\documentclass[leqno]{amsart}

\usepackage[foot]{amsaddr}

\usepackage{color}
\usepackage{geometry}
\usepackage{subfiles}
\usepackage{graphicx}
\usepackage[section]{placeins}
\usepackage{hyperref}
\usepackage{stmaryrd}
\usepackage{amsmath,amssymb,amsfonts,amsthm,textcomp}
\usepackage{mathtools}
\usepackage{tensor}
\usepackage{multicol}
\usepackage{subfig}
\usepackage{csquotes}
\usepackage{stackrel}
\usepackage{tikz}
\usepackage{pdfpages}
\graphicspath{
              {./Images/}
             }

\newfont{\tenbbb}{msbm10}
\newfont{\svnbbb}{msbm8}

\DeclareMathOperator*{\argmin}{arg\,min}
\DeclareMathOperator*{\argmax}{arg\,max}

\newcommand{\bs}[1]{\boldsymbol{#1}}
\newcommand{\cl}[1]{\mathcal{#1}}

\newcommand{\bb}[1]{\mathbb{#1}}

\newcommand{\red}[1]{{\color{red} #1}}

\newcommand{\di}[1]{\,\mathrm{d}{#1}}

\newcommand{\trans}{\scriptscriptstyle\mskip-1mu\top\mskip-2mu}

\newtheorem{rmk}{Remark}

\begin{document}

\title[Bayesian optimal experimental design with nuisance uncertainty]{Laplace-based strategies for Bayesian optimal experimental design with nuisance uncertainty}
\author{Arved Bartuska$^1$, Luis Espath$^2$, \& Ra\'{u}l Tempone$^{1,3,4}$}
\address{$^1$ Department of Mathematics, RWTH Aachen University, Geb\"{a}ude-1953 1.OG, Pontdriesch 14-16, 161, 52062 Aachen, Germany}
\address{$^2$ School of Mathematical Sciences, University of Nottingham, Nottingham, NG7 2RD, United Kingdom}
\address{$^3$ King Abdullah University of Science \& Technology (KAUST), Computer, Electrical and Mathematical Sciences \& Engineering Division (CEMSE), Thuwal 23955-6900, Saudi Arabia}
\address{$^4$ Alexander von Humboldt Professor in Mathematics for Uncertainty Quantification, RWTH Aachen University, Germany}
\email{arved.bartuska@gmail.com,espath@gmail.com,rtempone@gmail.com}

\date{\today}
\begin{abstract}
\noindent
%auto-ignore
Finding the optimal design of experiments in the Bayesian setting typically requires estimation and optimization of the expected information gain functional. This functional consists of one outer and one inner integral, separated by the logarithm function applied to the inner integral. When the mathematical model of the experiment contains uncertainty about the parameters of interest and nuisance uncertainty, (i.e., uncertainty about parameters that affect the model but are not themselves of interest to the experimenter), two inner integrals must be estimated. Thus, the already considerable computational effort required to determine good approximations of the expected information gain is increased further. The Laplace approximation has been applied successfully in the context of experimental design in various ways, and we propose two novel estimators featuring the Laplace approximation to alleviate the computational burden of both inner integrals considerably. The first estimator applies Laplace's method followed by a Laplace approximation, introducing a bias. The second estimator uses two Laplace approximations as importance sampling measures for Monte Carlo approximations of the inner integrals. Both estimators use Monte Carlo approximation for the remaining outer integral estimation. We provide four numerical examples demonstrating the applicability and effectiveness of our proposed estimators.
\\
\textbf{AMS subject classifications:}
$\cdot$
65K05 % Numerical analysis - Mathematical programming methods
$\cdot$
65C05 % Numerical analysis - Monte Carlo methods
$\cdot$

\end{abstract}

\maketitle

\tableofcontents                        % Print table of contents

%-------------------------------------------------------------------------------%

\section{Introduction}
The goal of optimal experimental design (OED) \cite{Cha95} is to provide designs for a given scientific experiment for which the expected information gain (EIG) \cite{Lin56} is optimal.
In the Bayesian setting, the knowledge of the experimenter about the parameters of interest before conducting the experiment is expressed as the prior probability density (pdf) of the parameters of interest, whereas the knowledge after the experiment is expressed as the posterior pdf, conditioned on the data. The increase in knowledge (or the reduction in uncertainty) is given by the Kullback--Leibler divergence \cite{Kul59, Kul51} between the posterior and prior.
The EIG can be expressed as the expected Kullback--Leibler divergence over the data. A higher EIG indicates that the data obtained from the experiment are expected to provide more information on the parameters of interest. Therefore, optimizing the EIG reduces the uncertainty about the parameters of interest.

If additional uncertainty is present in the model via parameters not of interest to the experimenter, we refer to them as nuisance parameters and marginalize them before optimizing the EIG. Thus, we find a design that only reduces the uncertainty about the parameters of interest, not the nuisance parameters. This approach ultimately maximizes the amount of information the experiment provides about the parameters of interest while keeping track of the uncertainty introduced by the nuisance parameters. A similar setting was considered in \cite{Ale22, Bar22, Fen19}. For a more general background on nuisance uncertainty, see \cite{Ber79, Lev14, Lie13, Pol88}.

There is often no closed-form expression for the EIG; thus, it must be estimated numerically \cite{Rya03}, which involves estimating nested integrals separated by the logarithm function and typically entails high computational costs. Marginalizing nuisance parameters introduces a second inner integral, increasing the computational cost further. The Laplace approximation \cite{Sti86, Tie86, Tie89, Kas90, Fri07} of the inner integral in combination with the Monte Carlo (MC) method is an effective tool for estimating the EIG \cite{Bec18, Lon13}. In this work, we extend this approach to the case with additional nuisance uncertainty and derive Laplace approximations for two distinct estimators. The first estimator uses Laplace's method followed by a Laplace approximation in a nested fashion to approximate the posterior pdf of the parameters of interest directly, and then uses the MC method for the outer integral. This method introduces a bias relative to the number of experiments performed, which is considered fixed from a modeling perspective. The second estimator uses two separate Laplace approximations as importance sampling measures for the two inner integrals, which are approximated using the MC method.

Our approach differs from that used in \cite{Fen19}, as the Laplace approximation is incorporated centrally in both estimators, whereas \cite{Fen19} developed a gradient-free importance sampling scheme to reduce variance in their proposed MC estimator. Similar to \cite{Fen19}, we consider reducible nuisance uncertainty, which is accessible to Bayesian updates. On the other hand, \cite{Ale22} considers irreducible nuisance uncertainty. The studies \cite{Ale22} and \cite{Bar22} assumed Gaussianity for the distribution of the nuisance parameters. Such assumptions are not present in \cite{Fen19} and the present work. The small-noise approximation derived in \cite{Bar22} also further restricts the size of the nuisance uncertainty to be considered.

In \cite{Eng22}, an approximate Laplace importance sampling for EIG estimation was proposed, where outer samples were used to compute the Hessian for the Laplace approximation rather than using the maximum a posteriori (MAP). A double-importance sampling scheme for EIG with nuisance uncertainty based on the Laplace approximation was first proposed in \cite{Eng22}. For the importance sampling density of the nuisance parameters conditioned on the data and the parameters of interest, the joint density of nuisance parameters and parameters of interest was first approximated by the standard Laplace approximation. The conditional density was then computed based on the formula for the conditional normal density. Instead, in this work we propose a novel Laplace approximation that is directly applied to the conditional density of the nuisance parameters. This approach requires additional numerical approximation methods, namely, the MAP and related Hessian for this conditional density; however, we also expect to obtain a more robust importance sampling density from our approach. Moreover, we demonstrate that the cost to obtain importance sampling densities does not affect the asymptotic cost of the overall estimator.

The work \cite{Ove17} considers normal-based Laplace approximations also for utility functions different from the EIG, called therein the expected Shannon information gain or expected self-information loss. Moreover, model uncertainty is also discussed in \cite{Ove17}, which can be interpreted as an application of nuisance parameters. However, the setting discussed in \cite{Ove17} is of a finite and countable number of possible models. Integral approximations, as required in the current manuscript for nuisance parameters modeled by a continuous random variable, are not discussed in their setting.

For asymptotic error bounds on the Laplace approximation, see \cite{Sch20, Wac17}. Pre-asymptotic error bounds and the effect of the nonlinearity of the experiment model on the data likelihood were investigated in \cite{Hel22}. For additional errors from estimating the mean and covariance in the Laplace approximation, see \cite{Spo22}.

The rest of this work is structured as follows. In Section~\ref{sec:Bayesian.formulation}, we introduce the Bayesian formulation and the standard double-loop MC (DLMC) estimator with nuisance uncertainty, which serves as a reference method for the novel estimators. Next, in Section~\ref{sec:Bias.stat}, we provide an error analysis, computational work analysis, and a method to obtain the required number of outer and inner samples to achieve a certain error tolerance for the DLMC estimator with nuisance uncertainty. As the main contribution of this work, we derive the Laplace approximations necessary for the MC double-Laplace (MC2LA) estimator in Section~\ref{sec:Laplace.approximation}. Subsequently, in Section~\ref{sec:Expected.information}, we develop the MC2LA estimator and a probabilistic error bound. Next, we develop the DLMC double importance sampling (DLMC2IS) estimator. Finally, in Section~\ref{sec:Numerical.results}, we demonstrate the effectiveness of the proposed estimators on four numerical examples. The first example provides insight into the effect of nuisance parameters on the optimal design of an experiment. An example from pharmacokinetics indicates the possibility for optimization in a larger design space of 15 dimensions. The third and fourth examples reveal the applicability of the estimators to electrical impedance tomography (EIT) experiments involving a finite element approximation of the underlying partial differential equation (PDE).

\section{Bayesian formulation}\label{sec:Bayesian.formulation}
Before we discuss the Bayesian setting, we must first describe the experiment in mathematical terms. For this purpose, we consider an additive data model:
\begin{equation}\label{eq:data.model}
\bs{y}_i(\bs{\xi})=\bs{g}(\bs{\xi},\bs{\theta}_t,\bs{\phi}_t)+\bs{\epsilon}_i,
\end{equation}
where
\begin{itemize}
  \item $\bs{y}_i\in\bb{R}^{d_y}$ is the observed data vector, $i=1,\hdots,N_e$;
  \item $\bs{Y}=(\bs{y}_1,\hdots,\bs{y}_{N_e}) \in \bb{R}^{d_y\times N_e}$ is the data for each experiment;
  \item $N_e$ is the number of observations for a specific experiment setup;
  \item $\bs{\theta}_t\in\bb{R}^{d_{\theta}}$ is the true vector of the parameters of interest;
  \item $\bs{\phi}_t\in\bb{R}^{d_{\phi}}$ is the true vector of the nuisance parameters;
	\item $\bs{\xi}\in\bb{R}^{d_{\xi}}$ is the vector of the design parameters;
  \item $\bs{g}\colon \bb{R}^{d_{\xi}}\times \bb{R}^{d_{\theta}}\times \bb{R}^{d_{\phi}} \to\bb{R}^{d_y}$ is the deterministic model of the experiment;
  \item $\bs{\epsilon}_i\in\bb{R}^{d_y}$ is the error vector assumed to be Gaussian $\bs{\epsilon}_{i}\stackrel{\mathrm{iid}}{\sim}\cl{N}(\bs{0},\bs{\Sigma}_{\bs{\epsilon}})$, $i=1,\hdots,N_e$,
\end{itemize}
where iid refers to independent and identically distributed.
The knowledge about the parameters of interest $\bs{\theta}$ before the experiment is encompassed by the prior pdf $\pi(\bs{\theta})$, whereas the knowledge after the experiment is given by the posterior pdf $\pi(\bs{\theta}|\bs{Y})$. We make no general assumptions on the independence of $\bs{\theta}$ and $\bs{\phi}$; therefore, we also consider the joint distribution $\pi(\bs{\theta}, \bs{\phi})$. To complete the Bayesian setup, we consider the likelihood $p(\bs{Y}|\bs{\theta})$ and evidence $p(\bs{Y})$. The latter two pdfs are denoted by $p$ rather than $\pi$ to distinguish between pdfs for model parameters and data. Moreover, the posterior, likelihood, and the evidence are conditioned on the design $\bs{\xi}$, but we omit this dependence for concision.
In this context, the Bayes formula for the posterior reads
\begin{align}\label{eq:posterior}
\pi(\bs{\theta}|\bs{Y}) ={}&
 \dfrac{\pi(\bs{\theta}) p(\bs{Y}|\bs{\theta})}{p(\bs{Y})},\nonumber\\
 ={}&\dfrac{\pi(\bs{\theta})\int_{\bs{\Phi}} p(\bs{Y}|\bs{\theta},\bs{\phi})\pi(\bs{\phi}|\bs{\theta})\di{}\bs{\phi}}{\int_{\bs{\Theta}}\int_{\bs{\Phi}}p(\bs{Y}|\bs{\theta},\bs{\phi})\pi(\bs{\theta},\bs{\phi})\di{}\bs{\phi}\di{}\bs{\theta}},
\end{align}
where the likelihood terms are obtained through marginalization. The likelihood conditioned on $\phi$ results from the data assumption \eqref{eq:data.model}, as
\begin{equation}\label{eq:likelihood}
p(\bs{Y}|\bs{\theta},\bs{\phi})\coloneqq\det(2\pi\bs{\Sigma}_{\bs{\varepsilon}})^{-\frac{N_e}{2}}\exp\left(-\frac{1}{2}\sum_{i=1}^{N_e} \bs{r}(\bs{y}_i,\bs{\theta},\bs{\phi}) \cdot \bs{\Sigma}_{\bs{\varepsilon}}^{-1}\bs{r}(\bs{y}_i,\bs{\theta},\bs{\phi})\right),
\end{equation}
where
\begin{equation}\label{eq:residual}
\bs{r}(\bs{y}_i,\bs{\theta},\bs{\phi})\coloneqq\bs{y}_i-\bs{g}(\bs{\theta},\bs{\phi})=\bs{g}(\bs{\theta}_t,\bs{\phi}_t)+\bs{\epsilon}_i-\bs{g}(\bs{\theta},\bs{\phi})
\end{equation}
is the data residual. The amount of information gained about the parameters of interest $\bs{\theta}$ from the experiment is expressed as the Kullback--Leibler divergence \cite{Kul59, Sha48} given by
\begin{equation}\label{Dkl}
  D_{\mathrm{KL}}=\int_{\bs{\Theta}}[\log(\pi(\bs{\theta}|\bs{Y}))-\log(\pi(\bs{\theta}))]\pi(\bs{\theta}|\bs{Y})\di{}\bs{\theta}.
\end{equation}

The goal of OED is to determine a design $\bs{\xi}$ for which the experiment provides informative data; therefore, we consider the expectation of \eqref{Dkl} over the data $\bs{Y}$, yielding the EIG

\begin{equation}\label{eq:EIG.post}
  I=\int_{\cl{Y}}\int_{\bs{\Theta}}[\log(\pi(\bs{\theta}|\bs{Y}))-\log(\pi(\bs{\theta}))]\pi(\bs{\theta}|\bs{Y})\di{}\bs{\theta} p(\bs{Y})\di{}\bs{Y}.
\end{equation}
This quantity only depends on the design $\bs{\xi}$, as all other dependencies have been marginalized. We rewrite this expression in terms of the known likelihood function \eqref{eq:likelihood} conditioned on $\phi$ using the Bayes formula \eqref{eq:posterior}
\begin{align}\label{eq:EIG.definition}
I={}&\int_{\bs{\Theta}}\int_{\cl{Y}}[\log(p(\bs{Y}|\bs{\theta}))-\log(p(\bs{Y}))]p(\bs{Y}|\bs{\theta}) \di{}\bs{Y}\pi(\bs{\theta})\di{}\bs{\theta},\nonumber\\
={}&\int_{\bs{\Theta}}\int_{\bs{\Phi}}\int_{\cl{Y}}\Bigg[\log\left(\int_{\bs{\Phi}} p(\bs{Y}|\bs{\theta},\bs{\varphi})\pi(\bs{\varphi}|\bs{\theta})\di{}\bs{\varphi}\right)\nonumber\\[4pt]
  {}&-\log\left(\int_{\bs{\Theta}}\int_{\bs{\Phi}}p(\bs{Y}|\bs{\vartheta},\bs{\varphi})\pi(\bs{\vartheta},\bs{\varphi})\di{}\bs{\varphi}\di{}\bs{\vartheta}\right)\Bigg]p(\bs{Y}|\bs{\theta},\bs{\phi})\di{}\bs{Y}\pi(\bs{\theta},\bs{\phi})\di{}\bs{\phi}\di{}\bs{\theta},
\end{align}
where $\bs\vartheta$ and $\bs\varphi$ are dummy variables to distinguish them from $\bs{\theta}$ and $\bs{\phi}$ from the outer integrals. As a reference to approximate expression \eqref{eq:EIG.definition}, we introduce the DLMC estimator with two inner loops:
\begin{equation}\label{EIG}
I_{\rm{DL}}= \frac{1}{N}\sum_{n=1}^N\left[\log\left(\frac{1}{M_1}\sum_{m=1}^{M_1}p(\bs{Y}^{(n)}|\bs{\theta}^{(n)},\bs{\varphi}^{(n,m)})\right)-\log\left(\frac{1}{M_2}\sum_{k=1}^{M_2}p(\bs{Y}^{(n)}|\bs{\vartheta}^{(n,k)},\bs{\varphi}^{(n,k)})\right)\right].
\end{equation}
The samples are drawn as follows. First, we sample $(\bs{\theta}^{(n)},\bs{\phi}^{(n)})\stackrel{\mathrm{iid}}{\sim}\pi(\bs{\theta},\bs{\phi})$. This allows us to sample $\bs{Y}^{(n)}\stackrel{\mathrm{iid}}{\sim}p(\bs{Y}|\bs{\theta}^{(n)},\bs{\phi}^{(n)})$, $1\leq n\leq N$. Next, we sample $\bs{\varphi}^{(n,m)}\stackrel{\mathrm{iid}}{\sim}\pi(\bs{\phi}|\bs{\theta}^{(n)})$, $1\leq m\leq M_1$, and $(\bs{\vartheta}^{(n,k)},\bs{\varphi}^{(n,k)})\stackrel{\mathrm{iid}}{\sim}\pi(\bs{\theta},\bs{\phi})$, $1\leq k\leq M_2$.

\section{Bias, statistical error, and optimal number of samples}\label{sec:Bias.stat}
To obtain the optimal number of samples for the outer and inner loops, we aim to minimize the total work to achieve a certain prescribed tolerance in the EIG estimator. Thus, we first stipulate that the total work $W$ of computing the estimator $I_{\rm{DL}}$ in expression \eqref{EIG} is given by
\begin{equation}\label{Work}
W(I_{\rm{DL}})\propto N(M_1+M_2).
\end{equation}
Moreover, we split the error of the estimator into the bias and statistical errors:
\begin{equation}
|I_{\rm{DL}}-I|\leq\underbrace{|\mathbb{E}[I_{\rm{DL}}]-I|}_{\text{bias error}}+\underbrace{|I_{\rm{DL}}-\mathbb{E}[I_{\rm{DL}}]|}_{\text{statistical error}}.
\end{equation}
Emulating the reckoning by \cite{Fen19}, with a second-order Taylor approximation, the bias component reads as
\begin{equation}\label{eq:bias.taylor}
|\mathbb{E}[I_{\rm{DL}}]-I|\approx \left|\frac{C_2}{M_2}-\frac{C_1}{M_1}\right|.
\end{equation}
In \eqref{eq:bias.taylor}, we have the difference between two positive terms; thus, these terms could cancel each other. In this scenario, the error estimation depends on the ignored higher-order terms. Thus, we lose all knowledge of the error estimates. To avoid this, we use the triangular inequality to arrive at a slightly more conservative error estimate:
\begin{align}\label{eq:bias.estimates}
|\mathbb{E}[I_{\rm{DL}}]-I|\approx{}& \left|\frac{C_2}{M_2}-\frac{C_1}{M_1}\right|,\nonumber\\
\leq{}&\frac{C_2}{M_2}+\frac{C_1}{M_1}.
\end{align}

For the statistical error, we also recall \cite{Fen19} and express the variance of the EIG estimator as follows:
\begin{equation}\label{eq:variance.estimates}
\mathbb{V}[I_{\rm{DL}}]\approx \frac{D_3}{N}+\frac{D_1}{NM_1}+\frac{D_2}{NM_2},
\end{equation}
where the constants $C_i>0$ and $D_i>0$ depend on the design parameter $\bs\xi$. Although we base this part of the work on \cite{Fen19}, the error estimates differ because we use \eqref{eq:bias.estimates} instead of \eqref{eq:bias.taylor}. Furthermore, \cite{Fen19} consider that $M_1=M_2$ for simplicity, whereas we derive an allocation of inner samples that is proportional to the constants $C_1$ and $C_2$, following \cite[Appendix A]{Bar22}.

The constants appearing in \eqref{eq:bias.estimates} and \eqref{eq:variance.estimates} are given as follows (see Appendix~\ref{ap:bias}):
\begin{equation}\label{eq:C1}
C_1 = \frac{1}{2}\mathbb{E}\left[\frac{\mathbb{V}[p(\bs{Y}|\bs{\theta},\bs{\varphi})|\bs{Y},\bs{\theta}]}{p^2(\bs{Y}|\bs{\theta})}\right],
\end{equation}
\begin{equation}\label{eq:C2}
C_2 = \frac{1}{2}\mathbb{E}\left[\frac{\mathbb{V}[p(\bs{Y}|\bs{\vartheta},\bs{\varphi})|\bs{Y}]}{p^2(\bs{Y})}\right],
\end{equation}
\begin{equation}\label{eq:D1}
  D_1 = \mathbb{E}\left[\frac{\mathbb{V}[p(\bs{Y}|\bs{\theta},\bs{\varphi})|\bs{Y},\bs{\theta}]}{p^2(\bs{Y}|\bs{\theta})}\right],
\end{equation}
\begin{equation}\label{eq:D2}
D_2 = \mathbb{E}\left[\frac{\mathbb{V}[p(\bs{Y}|\bs{\vartheta},\bs{\varphi})|\bs{Y}]}{p^2(\bs{Y})}\right],
\end{equation}
and
\begin{equation}\label{eq:D3}
D_3 = \mathbb{V}\left[\log\left(\frac{p(\bs{Y}|\bs{\theta})}{p(\bs{Y})}\right)\right].
\end{equation}

For a given error tolerance $TOL>0$, we can introduce a splitting parameter $\kappa\in(0,1)$ to distribute the error between the bias and statistical error. Thus, aided by the central limit theorem, we arrive at the following constraints:
\begin{equation}\label{st1}
|\mathbb{E}[I_{\rm{DL}}]-I|\leq (1-\kappa)TOL,
\end{equation}
\begin{equation}\label{st2}
\mathbb{V}[I_{\rm{DL}}]\leq \left(\frac{\kappa TOL}{C_\alpha}\right)^2,
\end{equation}
where $C_\alpha=\Phi^{-1}(1-\alpha/2)$ is the constant depending on the chosen confidence level $\alpha$, and $\Phi^{-1}$ is the inverse cumulative distribution function of a standard normal random variable.

To determine the optimal number of samples, we formulate the following subproblem. We find $N$, $M_1$, $M_2$, and $\kappa$ arising from the minimization of the total work \eqref{Work} subject to the bias and statistical error constraints \eqref{st1} and \eqref{st2}:
\begin{equation}\label{eq:subproblem}
\begin{aligned}
N^\ast, M_1^\ast, M_2^\ast, \kappa^\ast \coloneqq {}& \argmin_{N, M_1, M_2, \kappa} N(M_1+M_2)\\[4pt]
\textrm{s.t.} \quad {}& \frac{C_1}{M_1}+\frac{C_2}{M_2} \le (1-\kappa)TOL,\\[4pt]
\quad {}& \frac{1}{N}\left(D_3+\frac{D_1}{M_1}+\frac{D_2}{M_2}\right) \le \left(\frac{\kappa TOL}{C_\alpha}\right)^2,\\[4pt]
{}& N, M_1, M_2 > 0, \\[4pt]
{}& 1 > \kappa > 0.
\end{aligned}
\end{equation}

In Appendix~\ref{ap:optimal.setting}, we solve this problem analytically. We present the optimal sample choices as follows:
\begin{equation}\label{eq:opt.samples}
\left\{
\begin{aligned}
{}& \kappa^\ast=\frac{8TOL+3D_3-\sqrt{16TOLD_3+9D_3^2}}{8TOL},\\[4pt]
{}& N^\ast=\frac{C_\alpha^2\left(D_3+2(1-\kappa^\ast)TOL\right)}{(\kappa^\ast)^2 TOL^2},\\[4pt]
{}& M_1^\ast=\frac{C_1+\sqrt{C_1C_2}}{(1-\kappa^\ast)TOL},\\[4pt]
{}& M_2^\ast=\frac{C_2+\sqrt{C_1C_2}}{(1-\kappa^\ast)TOL},
\end{aligned}
\right.
\end{equation}
whereas the total optimal work reads as
\begin{equation}
W^\ast\propto\left(\frac{C_\alpha^2\left(D_3+2(1-\kappa^\ast)TOL\right)}{(\kappa^\ast)^2 TOL^2}\right)\left(\frac{(\sqrt{C_1}+\sqrt{C_2})^2}{(1-\kappa^\ast)TOL}\right)\propto TOL^{-3}.
\end{equation}

As $TOL$ approaches infinity, $\kappa$ approaches 1. Thus, for large tolerances, only the statistical error is relevant. As $TOL$ approaches 0, $\kappa$ approaches 2/3, which can be observed from the fact that the statistical error \eqref{eq:variance.estimates} can be rewritten to include the bias error using \eqref{eq:C1} to \eqref{eq:D2} as follows:
\begin{equation}\label{eq:var.bias}
  \bb{V}[I_{\rm{DL}}]\approx \frac{1}{N}\left(D_3+2|\mathbb{E}[I_{\rm{DL}}]-I|\right).
\end{equation}
Thus, the statistical error always takes priority over the bias error.

When the bias error due to the numerical discretization of the forward problem is also considered, we obtain the following term for the average computational work:
 \begin{equation}\label{Work.bias}
 W(I_{\rm{DL}}^h)\propto N(M_1+M_2)h^{-\gamma},
 \end{equation}
where $h^{-\gamma}$ is proportional to the average work of evaluating the forward model $\bs{g}$ with discretization parameter $h$. The subproblem \eqref{eq:subproblem} becomes
\begin{equation}\label{eq:subproblem.h}
\begin{aligned}
N^\ast, M_1^\ast, M_2^\ast, h^\ast, \kappa^\ast \coloneqq {}& \argmin_{N, M_1, M_2, h, \kappa} N(M_1+M_2)h^{-\gamma}\\[4pt]
\textrm{s.t.} \quad {}& C_3h^\eta+ \frac{C_1}{M_1}+\frac{C_2}{M_2} \le (1-\kappa)TOL,\\[4pt]
\quad {}& \frac{1}{N}\left(D_3+\frac{D_1}{M_1}+\frac{D_2}{M_2}\right) \le \left(\frac{\kappa TOL}{C_\alpha}\right)^2,\\[4pt]
{}& N, M_1, M_2, h > 0, \\[4pt]
{}& 1 > \kappa > 0,
\end{aligned}
\end{equation}
where $\eta$ is the weak convergence rate of the discretization method. We obtain the solution
\begin{equation}
\left\{
\begin{aligned}
{}& \kappa^\ast=\eta\frac{8TOL\eta+3D_3\eta+D_3\gamma+4TOL\gamma-\sqrt{D_3\left(9D_3\eta^2+6D_3\eta\gamma+D_3\gamma^2+16\eta^2 TOL+8TOL\eta\gamma\right)}}{2TOL\left(4\eta^2+4\eta\gamma+\gamma^2\right)},\\[4pt]
{}& N^{\ast}=\frac{C_{\alpha}^2}{\kappa^2TOL}\left(\frac{D_3}{TOL}+2\left(1-\kappa\left(1+\frac{\gamma}{2\eta}\right)\right)\right),\\[4pt]
{}& M_1^{\ast}=\frac{C_1\left(1+\sqrt{\frac{C_2}{C_1}}\right)}{\left(1-\kappa\left(1+\frac{\gamma}{2\eta}\right)\right)TOL},\\[4pt]
{}& M_2^{\ast}=\frac{C_2\left(1+\sqrt{\frac{C_1}{C_2}}\right)}{\left(1-\kappa\left(1+\frac{\gamma}{2\eta}\right)\right)TOL},\\[4pt]
{}& h^{\ast}=\left(\frac{\gamma\kappa TOL}{2\eta C_3}\right)^{\frac{1}{\eta}}.
\end{aligned}
\right.
\end{equation}
For the derivation of this result, see Appendix~\ref{ap:optimal.setting.bias} and \cite{Bar22, Bec18}.

\section{Laplace approximation}\label{sec:Laplace.approximation}
\subsection{Laplace's method}\label{sec:Laplace.method}
The DLMC estimator is expensive and often suffers from numerical underflow \cite{Bec18}. Some inexpensive and accurate estimators have been proposed for the EIG using importance sampling \cite{Bec18} or multilevel techniques combined with importance sampling \cite{Bec20} to reduce the computational cost. In addition, the inner integral may be estimated directly using the Laplace approximation. In this work, we enhance the computational performance of the DLMC estimator with nuisance uncertainty by deriving suitable Laplace approximations to marginalize the likelihoods arising in the Bayes update and integrate the inner loop. The Laplace approximation for the Bayes update in this setting is novel. Through what follows, we demonstrate that the Laplace approximation for the inner loop differs from that in \cite{Lon13,Bec18}.

Next, using \eqref{eq:residual}, the posterior \eqref{eq:posterior} may be expressed as follows:
\begin{align}\label{eq:posterior.to.laplace}
\pi(\bs{\theta}|\bs Y)={}& \frac{\pi(\bs{\theta})}{p(\bs{Y})\det(2\pi\bs{\Sigma}_{\bs{\varepsilon}})^{\frac{N_e}{2}}}\int_{\bs\Phi}\exp\left(-\frac{1}{2}\sum_{i=1}^{N_e}\bs r(\bs{y}_i,\bs{\theta},\bs{\phi}) \cdot \bs{\Sigma}_{\bs{\varepsilon}}^{-1}\bs r(\bs{y}_i,\bs{\theta},\bs{\phi})\right)\pi(\bs{\phi}|\bs{\theta})\text{d}\bs{\phi},\nonumber\\
={}& \frac{\pi(\bs{\theta})}{p(\bs{Y})\det(2\pi\bs{\Sigma}_{\bs{\varepsilon}})^{\frac{N_e}{2}}}\int_{\bs\Phi}\exp\underbrace{\left(-\frac{1}{2}\sum_{i=1}^{N_e}\bs r(\bs{y}_i,\bs{\theta},\bs{\phi}) \cdot \bs{\Sigma}_{\bs{\varepsilon}}^{-1}\bs r(\bs{y}_i,\bs{\theta},\bs{\phi})+\log(\pi(\bs{\phi}|\bs{\theta}))\right)}_{:=-f(\bs\theta,\bs{\phi})}\text{d}\bs{\phi},\nonumber\\
={}& \frac{\pi(\bs{\theta})}{p(\bs{Y})\det(2\pi\bs{\Sigma}_{\bs{\varepsilon}})^{\frac{N_e}{2}}}\int_{\bs\Phi} e^{-f(\bs\theta,\bs{\phi})}\text{d}\bs{\phi}.
\end{align}
We approximate the integral in \eqref{eq:posterior.to.laplace} using Laplace's method. Assuming that $f(\bs{\theta},\bs{\phi})$ has a unique minimum in $\bs{\phi}$, and that its Hessian in $\bs{\phi}$ is negative definite for almost all $\bs{\theta}$, we write the second-order Taylor approximation of $f(\bs\theta,\bs{\phi})$ in $\bs\phi$ around its minimum in $\bs\phi$, denoted as $\bs{\hat{\phi}}(\bs{\theta})$ and given by
\begin{align}
\bs{\hat{\phi}}(\bs{\theta})={}&\argmin_{\bs{\phi}}f(\bs\theta,\bs{\phi})\nonumber\\
={}&\argmin_{\bs{\phi}}\left(\frac{1}{2}\sum_{i=1}^{N_e}\bs r(\bs{y}_i,\bs{\theta},\bs{\phi}) \cdot \bs{\Sigma}_{\bs{\varepsilon}}^{-1}\bs r(\bs{y}_i,\bs{\theta},\bs{\phi})-\log(\pi(\bs{\phi}|\bs{\theta}))\right),\nonumber\\
={}&\argmin_{\bs{\phi}}\Bigg(\frac{1}{2}N_e(\bs g(\bs{\theta}_t,\bs{\phi}_t)-\bs g(\bs{\theta},\bs{\phi})) \cdot \bs{\Sigma}_{\bs{\varepsilon}}^{-1}(\bs g(\bs{\theta}_t,\bs{\phi}_t)-\bs g(\bs{\theta},\bs{\phi})),\nonumber\\
{}&+\sum_{i=1}^{N_e}\bs\epsilon_i\cdot{}\bs{\Sigma}_{\bs{\varepsilon}}^{-1}(\bs g(\bs{\theta}_t,\bs{\phi}_t)-\bs g(\bs{\theta},\bs{\phi}))-\log(\pi(\bs{\phi}|\bs{\theta}))\Bigg).
\end{align}
The Taylor expansion in $\bs{\phi}$ reads as
\begin{equation}
  \tilde f(\bs\theta,\bs\phi)=\underbrace{f(\bs\theta,\bs{\hat{\phi}}(\bs{\theta}))}_{\text{const. in } \bs\phi}+\underbrace{\nabla_{\bs\phi}f(\bs\theta,\bs{\hat{\phi}}(\bs{\theta}))}_{=\bs{0}}\cdot(\bs\phi-\bs{\hat{\phi}}(\bs{\theta}))+\frac{1}{2}(\bs\phi-\bs{\hat{\phi}}(\bs{\theta}))\cdot\nabla_{\bs\phi}\nabla_{\bs\phi}f(\bs\theta,\bs{\hat{\phi}}(\bs{\theta}))(\bs\phi-\bs{\hat{\phi}}(\bs{\theta})),
\end{equation}
with the following terms:
\begin{equation}
f(\bs{\theta},\bs{\hat{\phi}}(\bs{\theta}))=\frac{1}{2}\sum_{i=1}^{N_e}\bs r(\bs{y}_i,\bs{\theta},\bs{\hat{\phi}}(\bs{\theta})) \cdot \bs{\Sigma}_{\bs{\varepsilon}}^{-1}\bs r(\bs{y}_i,\bs{\theta},\bs{\hat{\phi}}(\bs{\theta}))-\log(\pi(\bs{\hat{\phi}}(\bs{\theta})|\bs{\theta})),
\end{equation}
\begin{equation}
\nabla_{\bs\phi}f(\bs{\theta},\bs{\hat{\phi}}(\bs{\theta}))=-\sum_{i=1}^{N_e}\nabla_{\bs{\phi}} \bs{g}(\bs{\theta},\bs{\hat{\phi}}(\bs{\theta}))^{\trans} \bs{\Sigma}_{\bs{\varepsilon}}^{-1}\bs r(\bs{y}_i,\bs{\theta},\bs{\hat{\phi}}(\bs{\theta}))-\nabla_{\bs{\phi}}\log(\pi(\bs{\hat{\phi}}(\bs{\theta})|\bs{\theta})),
\end{equation}
and
\begin{align}
\nabla_{\bs\phi}\nabla_{\bs\phi}f(\bs{\theta},\bs{\hat{\phi}}(\bs{\theta}))={}&-\sum_{i=1}^{N_e}\nabla_{\bs{\phi}}\nabla_{\bs{\phi}} \bs{g}(\bs{\theta},\bs{\hat{\phi}}(\bs{\theta}))^{\trans}\bs{\Sigma}_{\bs{\varepsilon}}^{-1}\bs r(\bs{y}_i,\bs{\theta},\bs{\hat{\phi}}(\bs{\theta}))+N_e\nabla_{\bs{\phi}} \bs{g}(\bs{\theta},\bs{\hat{\phi}}(\bs{\theta}))^{\trans} \bs{\Sigma}_{\bs{\varepsilon}}^{-1}\nabla_{\bs{\phi}} \bs{g}(\bs{\theta},\bs{\hat{\phi}}(\bs{\theta}))\nonumber\\
{}&-\nabla_{\bs{\phi}}\nabla_{\bs{\phi}}\log(\pi(\bs{\hat{\phi}}(\bs{\theta})|\bs{\theta})),
\end{align}
where $\nabla_{\bs{\phi}} \bs{g}(\bs{\theta},\bs{\hat{\phi}}(\bs{\theta}))\in\bb{R}^{d_y\times d_{\phi}}$ is the Jacobian of $\bs{g}$ with respect to $\bs{\phi}$, evaluated at $\bs{\hat{\phi}}(\bs{\theta})$. In addition, $\nabla_{\bs{\phi}}\nabla_{\bs{\phi}} \bs{g}(\bs{\theta},\bs{\hat{\phi}}(\bs{\theta})))\in\bb{R}^{d_y\times d_{\phi}\times d_{\phi}}$ is the Hessian of $\bs{g}$ with respect to $\bs{\phi}$, evaluated at $\bs{\hat{\phi}}(\bs{\theta})$. Because $\bs{\hat{\phi}}(\bs{\theta})$ is the minimizer of $f(\bs\theta,\bs{\phi})$, $\nabla_{\bs\phi}f(\bs\theta,\bs{\hat{\phi}}(\bs{\theta}))=\bs{0}$.

In index notation, using Einstein's convention for summation with Greek letters, we obtain
\begin{equation}
f=\frac{1}{2} \sum_{i=1}^{N_e} \bs{r}_{\alpha} (\bs{\Sigma}_{\bs{\varepsilon}}^{-1})_{\alpha\beta} \bs{r}_{\beta} - \log(\pi(\bs{\hat{\phi}}(\bs{\theta})|\bs{\theta})),
\end{equation}
\begin{equation}
(\nabla_{\bs{\phi}}f)_{\gamma}=-\sum_{i=1}^{N_e}(\partial_\gamma \bs{g}_{\alpha}) (\bs{\Sigma}_{\bs{\varepsilon}}^{-1})_{\alpha\beta} \bs{r}_{\beta} - \partial_\gamma (\log(\pi(\bs{\hat{\phi}}(\bs{\theta})|\bs{\theta}))),
\end{equation}
and
\begin{equation}
(\nabla_{\bs{\phi}}\nabla_{\bs{\phi}}f)_{\gamma\zeta}=-\sum_{i=1}^{N_e}(\partial_\zeta \partial_\gamma \bs{g}_{\alpha}) (\bs{\Sigma}_{\bs{\varepsilon}}^{-1})_{\alpha\beta} \bs{r}_{\beta} + N_e(\partial_\gamma \bs{g}_{\alpha}) (\bs{\Sigma}_{\bs{\varepsilon}}^{-1})_{\alpha\beta} (\partial_\zeta \bs{g}_{\beta}) - \partial_\zeta \partial_\gamma (\log(\pi(\bs{\hat{\phi}}(\bs{\theta})|\bs{\theta}))).
\end{equation}

Laplace's method reads as
\begin{equation}\label{eq:first.laplace}
\int_{\bs{\Phi}} e^{-f(\bs{\theta},\bs{\phi})}\di{}\bs{\phi}\approx\frac{(2\pi)^{\frac{d_{\phi}}{2}}}{\det(\nabla_{\bs\phi}\nabla_{\bs\phi}f(\bs{\theta},\bs{\hat{\phi}}(\bs{\theta})))^{\frac{1}{2}}}e^{-f(\bs{\theta},\bs{\hat{\phi}}(\bs{\theta}))},
\end{equation}
yielding the following approximation of the posterior:
\begin{equation}\label{eq:laplace.marginalization}
\pi_{\bs{\hat{\phi}}}(\bs{\theta}|\bs Y) \coloneqq \frac{\pi(\bs{\theta})\pi(\bs{\hat{\phi}}(\bs{\theta})|\bs{\theta})(2\pi)^{\frac{d_{\phi}}{2}}}{p(\bs{Y})\det(2\pi\bs{\Sigma}_{\bs{\varepsilon}})^{\frac{N_e}{2}}\det(\nabla_{\bs\phi}\nabla_{\bs\phi}f(\bs\theta,\bs{\hat{\phi}}(\bs{\theta})))^{\frac{1}{2}}}\exp\left(-\frac{1}{2}\sum_{i=1}^{N_e}\bs r(\bs{y}_i,\bs{\theta},\bs{\hat{\phi}}(\bs{\theta})) \cdot \bs{\Sigma}_{\bs{\varepsilon}}^{-1}\bs r(\bs{y}_i,\bs{\theta},\bs{\hat{\phi}}(\bs{\theta}))\right).
\end{equation}
Moreover, \eqref{eq:laplace.marginalization} is not Gaussian unless the prior is Gaussian, and the model $\bs{g}$ is linear on $\bs{\theta}$.

\subsection{Laplace approximation}\label{sec:Laplace.approx}
Assuming that \eqref{eq:laplace.marginalization} has a unique minimum and negative definite Hessian as a function of $\bs{\theta}$, we follow the approach in \cite{Lon13} to obtain a Gaussian approximation of the posterior distribution, taking the negative logarithm of \eqref{eq:laplace.marginalization} and evolving it around the maximum a posteriori estimate $\bs{\hat{\theta}}$ up to second order:
\begin{align}\label{eq:log.likelihood}
F(\bs{\theta})\coloneqq{}&-\log(\pi_{\bs{\hat{\phi}}}(\bs{\theta}|\bs Y)),\nonumber\\
={}&\frac{1}{2}\sum_{i=1}^{N_e}\bs r(\bs{y}_i,\bs{\theta},\bs{\hat{\phi}}(\bs{\theta})) \cdot \bs{\Sigma}_{\bs{\varepsilon}}^{-1}\bs r(\bs{y}_i,\bs{\theta},\bs{\hat{\phi}}(\bs{\theta}))-h(\bs{\theta})+k(\bs{\theta})-\ell(\bs{\theta})+C,
\end{align}
where
\begin{equation*}
h(\bs{\theta})\coloneqq\log(\pi(\bs{\theta})),
\end{equation*}
\begin{equation*}
k(\bs{\theta})\coloneqq\frac{1}{2}\log(\det(\nabla_{\bs\phi}\nabla_{\bs\phi}f(\bs\theta,\bs{\hat{\phi}}(\bs{\theta})))),
\end{equation*}
and
\begin{equation*}
\ell(\bs{\theta})\coloneqq\log(\pi(\bs{\hat{\phi}}(\bs{\theta})|\bs{\theta})).
\end{equation*}
The last two terms are new compared to their approach, and $C$ is a constant.

The maximum a posteriori estimate $\bs{\hat{\theta}}$ of \eqref{eq:laplace.marginalization} is given by
\begin{align}\label{eq:theta.hat}
\bs{\hat{\theta}}={}&\argmax_{\bs{\theta}}\pi_{\bs{\hat{\phi}}}(\bs{\theta}|\bs Y),\nonumber\\
={}&\argmin_{\bs{\theta}}F(\bs{\theta}),\nonumber\\
={}&\argmin_{\bs{\theta}}\left[\frac{1}{2}\sum_{i=1}^{N_e}\bs r(\bs{y}_i,\bs{\theta},\bs{\hat{\phi}}(\bs{\theta})) \cdot \bs{\Sigma}_{\bs{\varepsilon}}^{-1}\bs r(\bs{y}_i,\bs{\theta},\bs{\hat{\phi}}(\bs{\theta}))-h(\bs{\theta})+k(\bs{\theta})-\ell(\bs{\theta})\right].
\end{align}
The Taylor expansion is given by
\begin{equation}\label{eq:second.taylor}
  \tilde F(\bs\theta)=\underbrace{F(\bs{\hat{\theta}})}_{\text{const. in } \bs\theta}+\underbrace{\nabla_{\bs\theta} F(\bs{\hat{\theta}})}_{=\bs{0}}\cdot(\bs\theta-\bs{\hat{\theta}})+\frac{1}{2}(\bs\theta-\bs{\hat{\theta}})\cdot\nabla_{\bs\theta}\nabla_{\bs\theta} F(\bs{\hat{\theta}})(\bs\theta-\bs{\hat{\theta}}).
\end{equation}
From \eqref{eq:log.likelihood}, we obtain
\begin{equation}
  \nabla_{\bs\theta} F(\bs{\bs{\hat{\theta}}})=-\sum_{i=1}^{N_e}\left(\nabla_{\bs{z}}\bs{g}(\bs{z}(\bs{\hat{\theta}}))\nabla_{\bs{\theta}}\bs{z}(\bs{\hat{\theta}})\right)^{\trans} \bs{\Sigma}_{\bs{\varepsilon}}^{-1}\bs{r}(\bs{y}_i,\bs{z}(\bs{\hat{\theta}}))-\nabla_{\bs\theta}h(\bs{\bs{\hat{\theta}}})+\nabla_{\bs\theta}k(\bs{\bs{\hat{\theta}}})-\nabla_{\bs{\theta}}\ell(\bs{\hat{\theta}}),
\end{equation}
by the chain rule, where we write
\begin{equation}
  \bs{z}(\bs{\hat{\theta}})\coloneqq (\bs{\hat{\theta}},\bs{\hat{\phi}}(\bs{\hat{\theta}}))\in\bb{R}^{d_{\theta}+ d_{\phi}},
\end{equation}
and
\begin{align}\label{eq:Ht}
  \nabla_{\bs\theta}\nabla_{\bs\theta} F(\bs{\bs{\hat{\theta}}})={}&-\sum_{i=1}^{N_e}\left(\nabla_{\bs{\theta}}\bs{z}(\bs{\hat{\theta}})^{\trans}\nabla_{\bs{z}}\nabla_{\bs{z}}\bs{g}(\bs{z}(\bs{\hat{\theta}}))\nabla_{\bs{\theta}}\bs{z}(\bs{\hat{\theta}})+ \nabla_{\bs{z}}\bs{g}(\bs{z}(\bs{\hat{\theta}}))\nabla_{\bs{\theta}}\nabla_{\bs{\theta}}\bs{z}(\bs{\hat{\theta}}) \right)^{\trans}\bs{\Sigma}_{\bs{\varepsilon}}^{-1}\bs{r}(\bs{y}_i,\bs{z}(\bs{\hat{\theta}}))\nonumber\\
  {}&+N_e\left(\nabla_{\bs{z}}\bs{g}(\bs{z}(\bs{\hat{\theta}}))\nabla_{\bs{\theta}}\bs{z}(\bs{\hat{\theta}})\right)^{\trans}\bs{\Sigma}_{\bs{\varepsilon}}^{-1}\left(\nabla_{\bs{z}}\bs{g}(\bs{z}(\bs{\hat{\theta}}))\nabla_{\bs{\theta}}\bs{z}(\bs{\hat{\theta}})\right)-\nabla_{\bs{\theta}}\nabla_{\bs{\theta}}h(\bs{\bs{\hat{\theta}}})+\nabla_{\bs{\theta}}\nabla_{\bs{\theta}}k(\bs{\hat{\theta}})-\nabla_{\bs{\theta}}\nabla_{\bs{\theta}}\ell(\bs{\hat{\theta}}).
\end{align}
Given \eqref{eq:log.likelihood}, \eqref{eq:theta.hat}, \eqref{eq:second.taylor}, and \eqref{eq:Ht} we find the Gaussian approximation of the posterior at $\bs{\hat{\theta}}$:
\begin{equation}\label{ptilde}
\tilde{\pi}(\bs{\theta}|\bs{Y})\coloneqq\frac{1}{\det(2\pi\bs{\Sigma})^{\frac{1}{2}}}\exp\left(-\frac{(\bs{\theta}-\bs{\hat\theta})\cdot\bs{\Sigma}^{-1}(\bs{\theta}-\bs{\hat\theta})}{2}\right),
\end{equation}
where $\bs{\Sigma}^{-1}=\nabla_{\bs\theta}\nabla_{\bs\theta} F(\bs{\bs{\hat{\theta}}})$. The first term in \eqref{eq:Ht} is $\cl{O}_{\bb{P}}(\sqrt{N_e})$\footnote{The notation $X_M=\cl{O}_{\bb{P}}(a_M)$ for a sequence of random variables $X_M$ and constants $a_M$ is as follows. For any $\epsilon>0$, there exists a finite $K(\epsilon)>0$ and finite $M_0>0$ such that $\bb{P}(|X_M|>K(\epsilon)|a_M|)<\epsilon$ holds for all $M\geq M_0$.}, and the second term is $\cl{O}_{\bb{P}}(N_e)$ (see \cite{Lon13} and Appendix~\ref{ap:order.maximizer}). The last two terms are of order $\cl{O}_{\bb{P}}(1)$.

\begin{rmk}[Differences and additional expenses in accounting for nuisance uncertainty in the Laplace approximation of the posterior]
In this setting, we obtain terms relating to $\bs{z}(\bs{\hat{\theta}})$, $k$, and $\ell$, requiring additional Jacobian and Hessian evaluations to adequately account for nuisance uncertainty.
\end{rmk}

\begin{rmk}[Uniqueness of the minimum]
    The Laplace approximation and Laplace's method can be applied to functions without a unique minimum, see \cite{Bor11} and \cite{Lon22}. However, the examples considered in Section~\ref{sec:Numerical.results} displayed behavior consistent with the assumptions stated at the beginning of Section~\ref{sec:Laplace.method} and Section~\ref{sec:Laplace.approx}.
\end{rmk}

\section{Expected Information Gain Estimators}\label{sec:Expected.information}
\subsection{Double Laplace approximation: Monte Carlo double Laplace}
The first estimator to compute the EIG uses Laplace's method and a Laplace approximation. We begin by rewriting the log ratio between posterior and prior:
\begin{equation}
	\log\left(\frac{\pi(\bs{\theta}|\bs Y)}{\pi(\bs{\theta})}\right)=\underbrace{\log\left(\frac{\pi(\bs{\theta}|\bs Y)}{\tilde{\pi}(\bs{\theta}|\bs Y)}\right)}_{:=\epsilon_{La}}+\log\left(\frac{\tilde{\pi}(\bs{\theta}|\bs Y)}{\pi(\bs{\theta})}\right).
\end{equation}
Using the second Laplace approximation given in \eqref{ptilde}, we obtain
\begin{equation}
 	\log\left(\frac{\pi(\bs{\theta}|\bs Y)}{\pi(\bs{\theta})}\right)=\epsilon_{La}-\frac{1}{2}\log(\det(2\pi\bs{\Sigma}))-\left(\frac{(\bs{\theta}-\bs{\hat\theta}) \cdot \bs{\Sigma}^{-1}(\bs{\theta}-\bs{\hat\theta})}{2}\right)-\log(\pi(\bs{\theta})).
\end{equation}
This formulation is identical to that \cite{Lon13}, and the only difference in the formulation without nuisance uncertainty is encoded in the covariance matrix $\bs{\Sigma}$.

Next, we write the Kullback--Leibler divergence \eqref{Dkl} as follows:
\begin{align}
	D_{\mathrm{KL}}={}&\int_{\bs\Theta}\log\left(\frac{\pi(\bs{\theta}|\bs Y)}{\pi(\bs{\theta})}\right)\tilde{\pi}(\bs{\theta}|\bs Y)\di{}\bs\theta+\underbrace{\int_{\bs\Theta}\log\left(\frac{\pi(\bs{\theta}|\bs Y)}{\pi(\bs{\theta})}\right)(\pi(\bs{\theta}|\bs Y)-\tilde{\pi}(\bs{\theta}|\bs Y))\di{}\bs\theta}_{\coloneqq \epsilon_{int}},\nonumber\\
	={}&\int_{\bs\Theta}\epsilon_{La}\tilde{\pi}(\bs{\theta}|\bs Y)\di{}\bs\theta+\int_{\bs\Theta}\left[-\frac{1}{2}\log(\det(2\pi\bs{\Sigma}))-\left(\frac{(\bs{\theta}-\bs{\hat{\theta}}) \cdot \bs{\Sigma}^{-1}(\bs{\theta}-\bs{\hat\theta})}{2}\right)-\log(\pi(\bs{\theta}))\right]\tilde{\pi}(\bs{\theta}|\bs Y)\di{}\bs\theta+\epsilon_{int},\nonumber\\
	={}&-\frac{1}{2}\log(\det(2\pi\bs{\Sigma}))-\frac{d_{\theta}}{2}-\log(\pi(\bs{\hat{\theta}}))-\frac{\text{tr}(\bs\Sigma\nabla_{\bs\theta}\nabla_{\bs\theta}\log(\pi(\bs{\hat{\theta}})))}{2}+\cl{O}_{\bb{P}}\left(\frac{1}{N_e^2}\right),
\end{align}
where $\epsilon_{La}, \, \epsilon_{int}=\cl{O}_{\bb{P}}\left(\frac{1}{N_e^2}\right)$ (see \cite[Appendices A, B, and C]{Lon13}) and $\log(\pi(\bs{\theta}))=\log(\pi(\bs{\hat{\theta}}))+\frac{\text{tr}(\bs\Sigma\nabla_{\bs\theta}\nabla_{\bs\theta}\log(\pi(\bs{\hat{\theta}})))}{2}+\cl{O}_{\bb{P}}\left(\frac{1}{N_e^2}\right)$ (see \cite[Appendix A]{Lon13}).

Taking the expected value over all $\bs Y$ and using the law of total probability results in the following EIG:
\begin{multline}
	I=\int_{\bs{\Theta}}\int_{\bs{\Phi}}\int_{\cl{Y}}\left[-\frac{1}{2}\log(\det(2\pi\bs{\Sigma}))-\frac{d_{\theta}}{2}-\log(\pi(\bs{\hat{\theta}}))-\frac{\text{tr}(\bs{\Sigma}\nabla_{\bs{\theta}}\nabla_{\bs{\theta}}\log(\pi(\bs{\hat{\theta}})))}{2}\right]\\
	\times p(\bs{Y}|\bs{\theta}, \bs{\phi})\di{}\bs{Y}\pi(\bs{\theta},\bs{\phi})\di{}\bs{\phi}\di{}\bs{\theta}+\cl{O}_{\bb{P}}\left(\frac{1}{N_e^2}\right),
\end{multline}
whereas its sample-based version reads as
\begin{equation}
	\hat{I}=\dfrac{1}{N}\sum_{n=1}^N\left[-\frac{1}{2}\log(\det(2\pi\bs{\Sigma}(\bs{\hat{\theta}}(\bs{Y}^{(n)}))))-\frac{d_{\theta}}{2}-\log(\pi(\bs{\hat{\theta}}(\bs{Y}^{(n)})))-\frac{\text{tr}(\bs{\Sigma}(\bs{\hat{\theta}}(\bs{Y}^{(n)}))\nabla_{\bs{\theta}}\nabla_{\bs{\theta}}\log(\pi(\bs{\hat{\theta}}(\bs{Y}^{(n)}))))}{2}\right],
\end{equation}
where $\bs{\Sigma}$ is the inverse Hessian \eqref{eq:Ht} evaluated at $\bs{\hat{\theta}}$. Both $\bs{\Sigma}$ and  $\bs{\hat{\theta}}$ depend on $\bs{Y}^{(n)}\stackrel{\mathrm{iid}}{\sim} p(\bs{Y}|\bs{\theta}^{(n)},\bs{\phi}^{(n)})$, where $(\bs{\theta}^{(n)},\bs{\phi}^{(n)})\stackrel{\mathrm{iid}}{\sim}\pi(\bs{\theta},\bs{\phi})$, $1\leq n\leq N$, are sampled from the prior. The optimal setting for this estimator was derived in \cite{Bec18}.

\subsection{Double Laplace-based Importance Sampling: Double Loop Monte Carlo double Importance Sampling}
Rather than directly approximating the inner integrals, we can also use the Laplace approximation as a change of measure in the EIG \eqref{eq:EIG.definition}. This method is known as importance sampling and reduces the variances \eqref{eq:C1}, \eqref{eq:C2}, \eqref{eq:D1}, and \eqref{eq:D2} in the inner MC loops, reducing the total work \eqref{Work}. We aim to estimate the following:
\begin{align}\label{eq:EIG.IS}
  I={}&\int_{\bs{\Theta}}\int_{\bs{\Phi}}\int_{\cl{Y}}\Bigg[\log\left(\int_{\bs{\Phi}} p(\bs{Y}|\bs{\theta},\bs{\varphi})\frac{\pi(\bs{\varphi}|\bs{\theta})}{\tilde{\pi}(\bs{\varphi}|\bs{Y},\bs{\theta})}\tilde{\pi}(\bs{\varphi}|\bs{Y},\bs{\theta})\di{}\bs{\varphi}\right)\nonumber\\[4pt]
  {}&-\log\left(\int_{\bs{\Theta}}\int_{\bs{\Phi}}p(\bs{Y}|\bs{\vartheta},\bs{\varphi})\frac{\pi(\bs{\vartheta},\bs{\varphi})}{\tilde{\pi}(\bs{\vartheta},\bs{\varphi}|\bs{Y})}\tilde{\pi}(\bs{\vartheta},\bs{\varphi}|\bs{Y})\di{}\bs{\varphi}\di{}\bs{\vartheta}\right)\Bigg]p(\bs{Y}|\bs{\theta},\bs{\phi})\di{}\bs{Y}\pi(\bs{\theta},\bs{\phi})\di{}\bs{\phi}\di{}\bs{\theta},
\end{align}
where $\tilde{\pi}(\bs{\varphi}|\bs{Y},\bs{\theta})$ is the Laplace approximation of $\pi(\bs{\varphi}|\bs{Y},\bs{\theta})$ and $\tilde{\pi}(\bs{\vartheta},\bs{\varphi}|\bs{Y})$ is the Laplace approximation of $\pi(\bs{\vartheta},\bs{\varphi}|\bs{Y})$. These Laplace approximations are different from the one developed in the previous section. The data $\bs{Y}=(\bs{y}_1,\ldots,\bs{y}_{N_e})$ in \eqref{eq:EIG.IS} can be decomposed into $\bs{y}_i=\bs{g}(\bs{\theta},\bs{\phi})+\bs{\epsilon}_i$, $1\leq i\leq N_e$; therefore, we write the distribution of $\bs{\varphi}$ as follows:
\begin{equation}
	\pi(\bs{\varphi}|\bs{Y},\bs{\theta})\propto\pi(\bs{\varphi}|\bs{\theta})p(\bs{Y}|\bs{\theta},\bs{\varphi}),
\end{equation}
where
\begin{equation}
	p(\bs{Y}|\bs{\theta},\bs{\varphi})\propto\exp\left(-\frac{1}{2}\sum_{i=1}^{N_e} \bs{r}(\bs{y}_i,\bs{\theta},\bs{\varphi}) \cdot \bs{\Sigma}_{\bs{\varepsilon}}^{-1}\bs{r}(\bs{y}_i,\bs{\theta},\bs{\varphi})\right),
\end{equation}
and
\begin{align}
\bs{r}(\bs{y}_i,\bs{\theta},\bs{\varphi})={}&\bs{y}_i-\bs{g}(\bs{\theta},\bs{\varphi}),\nonumber\\
={}&\bs{g}(\bs{\theta},\bs{\phi})+\bs{\epsilon}_i-\bs{g}(\bs{\theta},\bs{\varphi}), \quad 1\leq i\leq N_e.
\end{align}
Following the steps in \cite{Lon13}, we arrive at
\begin{equation}
	\tilde{\pi}(\bs{\varphi}|\bs{Y},\bs{\theta})\coloneqq\frac{1}{\det(2\pi\bs{\Sigma}_{\bs{\phi}})^{\frac{1}{2}}}\exp\left(-\frac{(\bs{\varphi}-\bs{\hat{\varphi}})\cdot\bs{\Sigma}_{\bs{\phi}}^{-1}(\bs{\varphi}-\bs{\hat{\varphi}})}{2}\right),
\end{equation}
where
\begin{equation}
	\bs{\hat{\varphi}}=\argmin_{\bs{\varphi}}\left[\frac{1}{2}\sum_{i=1}^{N_e}\bs r(\bs{y}_i,\bs{\theta},\bs{\varphi}) \cdot \bs{\Sigma}_{\bs{\varepsilon}}^{-1}\bs r(\bs{y}_i,\bs{\theta},\bs{\varphi})-\log(\pi(\bs{\varphi}|\bs{\theta}))\right],
\end{equation}
and
\begin{equation}\label{eq:Sigma.phi}
	\bs{\Sigma}_{\bs{\phi}}^{-1}=-\sum_{i=1}^{N_e}\nabla_{\bs{\varphi}}\nabla_{\bs{\varphi}}\bs{g}(\bs{\theta},\bs{\hat{\varphi}})^{\trans}\bs{\Sigma}_{\bs{\varepsilon}}^{-1}\bs r(\bs{y}_i,\bs{\theta},\bs{\hat{\varphi}})+N_e\nabla_{\bs{\varphi}}\bs{g}(\bs{\theta},\bs{\hat{\varphi}})^{\trans} \bs{\Sigma}_{\bs{\varepsilon}}^{-1}\nabla_{\bs{\varphi}}\bs{g}(\bs{\theta},\bs{\hat{\varphi}})-\nabla_{\bs{\phi}}\nabla_{\bs{\phi}}\log(\pi(\bs{\hat{\varphi}}|\bs{\theta})).
\end{equation}
We drop the first term in \eqref{eq:Sigma.phi}, as it is $\cl{O}_{\bb{P}}\left(\sqrt{N_e}\right)$. As for the joint density $\pi(\bs{\vartheta},\bs{\varphi}|\bs{Y})$, we write $\bs z\coloneqq(\bs{\vartheta},\bs{\varphi})\in\bb{R}^{d_{\theta}+ d_{\phi}}$ to obtain the Laplace approximation:
\begin{equation}\label{eq:Laplace.z}
	\tilde{\pi}(\bs{z}|\bs{Y})\coloneqq\frac{1}{\det(2\pi\bs{\Sigma}_{\bs{z}})^{\frac{1}{2}}}\exp\left(-\frac{(\bs{z}-\bs{\hat{z}})\cdot\bs{\Sigma}_{\bs{z}}^{-1}(\bs{z}-\bs{\hat{z}})}{2}\right),
\end{equation}
with
\begin{equation}
	\bs{\hat{z}}=\argmin_{\bs{z}}\left[\frac{1}{2}\sum_{i=1}^{N_e}\bs r(\bs{y}_i,\bs{z}) \cdot \bs{\Sigma}_{\bs{\varepsilon}}^{-1}\bs r(\bs{y}_i,\bs{z})-\log(\pi(\bs{z}))\right],
\end{equation}
and
\begin{equation}\label{eq:Sigma.z}
	\bs{\Sigma}_{\bs{z}}^{-1}=-\sum_{i=1}^{N_e}\nabla_{\bs{z}}\nabla_{\bs{z}}\bs{g}(\bs{\hat{z}})^{\trans}\bs{\Sigma}_{\bs{\varepsilon}}^{-1}\bs r(\bs{y}_i,\bs{\hat{z}})+N_e\nabla_{\bs{z}}\bs{g}(\bs{\hat{z}})^{\trans} \bs{\Sigma}_{\bs{\varepsilon}}^{-1}\nabla_{\bs{z}}\bs{g}(\bs{\hat{z}})-\nabla_{\bs{z}}\nabla_{\bs{z}}\log(\pi(\bs{\hat{z}})).
\end{equation}
This last Laplace approximation \eqref{eq:Laplace.z} has the same shape as that derived in \cite{Lon13}. We again drop the first term in \eqref{eq:Sigma.z}, as it is $\cl{O}_{\bb{P}}\left(\sqrt{N_e}\right)$.
This method provides the DLMC2IS estimator, as follows:
\begin{align}
\hat{I}={}& \frac{1}{N}\sum_{n=1}^N\left[\log\left(\frac{1}{M_1}\sum_{m=1}^{M_1}p(\bs Y^{(n)}|\bs{\theta}^{(n)},\bs{\varphi}^{(n,m)})\frac{\pi(\bs{\varphi}^{(n,m)})}{\tilde{\pi}(\bs{\varphi}^{(n,m)}|\bs{Y}^{(n)},\bs{\theta}^{(n)})}\right)\right.\nonumber\\
{}&\left.-\log\left(\frac{1}{M_2}\sum_{k=1}^{M_2}p(\bs Y^{(n)}|\bs{\vartheta}^{(n,k)},\bs{\varphi}^{(n,k)})\frac{\pi(\bs{\vartheta}^{(n,k)},\bs{\varphi}^{(n,k)})}{\tilde{\pi}(\bs{\vartheta}^{(n,k)},\bs{\varphi}^{(n,k)}|\bs{Y}^{(n)})}\right)\vphantom{\log\left(\frac{1}{M_1}\sum_{m=1}^{M_1}p(\bs Y^{(n)}|\bs{\theta}^{(n)},\bs{\varphi}^{(n,m)})\frac{\pi(\bs{\varphi}^{(n,m)})}{\tilde{\pi}(\bs{\varphi}^{(n,m)}|\bs{Y}^{(n)},\bs{\theta}^{(n)})}\right)}\right].
\end{align}
First, we sample $(\bs{\theta}^{(n)},\bs{\phi}^{(n)})\stackrel{\mathrm{iid}}{\sim}\pi(\bs{\theta},\bs{\phi})$, then we sample $\bs Y^{(n)}\stackrel{\mathrm{iid}}{\sim}p(\bs Y|\bs{\theta}^{(n)},\bs{\phi}^{(n)})$, $\bs{\varphi}^{(n,m)}\stackrel{\mathrm{iid}}{\sim}\tilde{\pi}(\bs{\varphi}|\bs{Y}^{(n)},\bs{\theta}^{(n)})$ and  $(\bs{\vartheta}^{(n,k)},\bs{\varphi}^{(n,k)})\stackrel{\mathrm{iid}}{\sim}\tilde{\pi}(\bs{\vartheta},\bs{\varphi}|\bs{Y}^{(n)})$.
The derivation of the optimal setting for the DLMC2IS estimator follows the same basic structure as for the DLMC estimator. Laplace-based importance sampling helps reduce the variance of the inner MC approximation, and thus affects the constants appearing in the bias and variance bounds presented in Section~\ref{sec:Bias.stat}. Moreover, the total work for the DLMC2IS estimator includes the cost of numerically estimating the MAP and Hessian for each outer sample \cite{Bec18}.  
The asymptotic computational work of the DLMC estimator is
\begin{equation}
    W_{DLMC}=\cl{O}\left(N\times(M_1+M_2)\times h^{-\gamma}\right),
\end{equation}
where $N,M_1,M_2$, and $h$ are expressed as functions of the error tolerance $TOL$. Preasymptotically, this work may be expressed as
\begin{equation}
    W_{DLMC}=c_1N\times(c_{2}M_1+c_{3}M_2+1)\times c_4h^{-\gamma},
\end{equation}
for appropriately chosen constants $c_1,c_2,c_3,c_4$ independent of $TOL$. The factor 1 enters as the cost of evaluating the experiment model once per outer sample to generate the data $\bs{Y}$, which can be ignored asymptotically. For the DLMC2IS estimator, it holds that
\begin{equation}
    W_{DLIS}=c_1N\times(\tilde{c}_{2}M_1+\tilde{c}_{3}M_2+J)\times c_4h^{-\gamma},
\end{equation}
where typically $\tilde{c}_2\ll c_2$, $\tilde{c}_3\ll c_3$ due to importance sampling, and $J\geq 2$ represents the cost to numerically estimate MAP and Hessian for each outer sample. This factor $J$, as well as $\tilde{c}_1$ and $\tilde{c}_2$, are independent of $TOL$, thus, it follows that
\begin{align}
    c_1N\times(\tilde{c}_{2}M_1+\tilde{c}_{3}M_2+J)\times c_4h^{-\gamma}{}&=c_1N\times(\tilde{c}_{2}M_1+\tilde{c}_{3}M_2)\times c_4h^{-\gamma}+c_1N\times J\times c_4h^{-\gamma},\nonumber\\
    {}&\geq (1+\epsilon)\times c_1N\times(\tilde{c}_{2}M_1+\tilde{c}_{3}M_2)\times c_4h^{-\gamma},
\end{align}
for any $\epsilon>0$ as $TOL\to 0$, where the last term in the first line is of higher order in $TOL$. The asymptotic work of the estimator is therefore influenced by a reduction in a multiplicative error term proportional to the combined importance sampling effects. The cost for finding MAPs and Hessians only enters as an arbitrary $\epsilon>0$.
Numerical results in Figure~\ref{fig:ex1.evt} (Panel (A)) show that our estimator is conservative, implying the number of samples used was larger than necessary for a given tolerance.
\rmk{[Importance sampling for a prior distribution with compact support] If the support of the prior distribution is compact, using a Gaussian distribution for importance sampling can result in samples outside this domain (see \cite{Bis16}). However, this problem is rarely observed in cases where the covariance of the importance sampling distribution is highly concentrated.
}
\rmk{[Numerical estimation of MAPs and Hessians]
All MAPs and Hessians in the following section were estimated numerically. Specific experiment models allow for closed-form expressions or the use of automatic differentiation; however, this was beyond the scope of this work.
}

\section{Numerical Results}\label{sec:Numerical.results}
\subsection{Linear Gaussian Example}\label{ex:1}
\begin{figure}[ht]
		\includegraphics[width=0.6\textwidth]{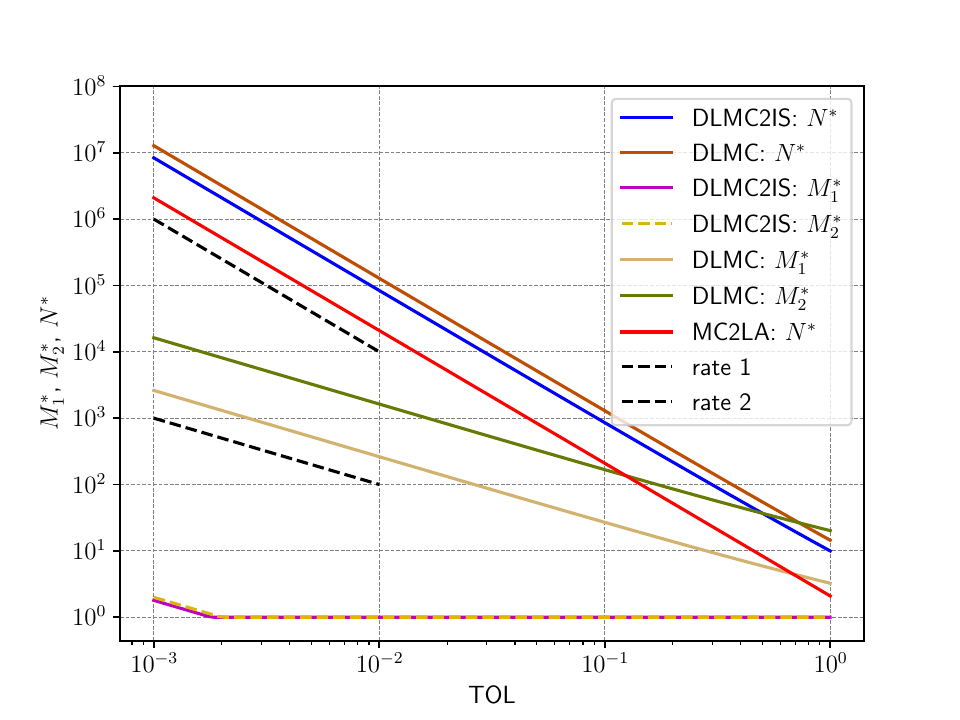}
	\caption{Example 1: Optimal number of outer ($N^{\ast}$) and inner ($M_1^{\ast}$, $M_2^{\ast}$) samples vs. tolerance $TOL$ for the DLMC, DLMC2IS, and MC2LA estimators.}
	\label{fig:ex1.opt}
\end{figure}

This example is based on a similar formulation used by \cite{Fen19} to demonstrate the effects of nuisance parameters on the optimal design $\bs{\xi}$. We assume a linear model,
\begin{equation}
  \bs{y}(\xi) = \bs{g}(\xi,\theta,\phi)+\bs{\epsilon},
\end{equation}
where
\begin{align}\label{eq:ex1.model}
  \bs{g}(\xi,\theta,\phi)={}&\begin{pmatrix}
  	\xi& 0 \\ 0&(1-\xi)
  \end{pmatrix}\begin{pmatrix}
  	\theta \\ \phi
  \end{pmatrix},\nonumber\\
  ={}&\begin{pmatrix}
    \xi\theta\\(1-\xi)\phi
  \end{pmatrix},
\end{align}
and the parameters are sampled from the following distributions:
\begin{equation}
  \theta\sim\cl{N}(0,1),\quad \phi\sim\cl{N}(0,10^{-2}),\quad \bs{\epsilon}\sim\cl{N}(\bs{0},\bs{\Sigma}_{\bs{\epsilon}}),\quad \bs{\Sigma}_{\bs{\epsilon}}=\begin{pmatrix}
    10^{-2}&0\\0&10^{-2}
  \end{pmatrix}.
\end{equation}
The vector $(0, 0)^{\trans}$ is denoted by $\bs{0}$. The design $\xi$ is chosen from the interval $(0,1)$, and we only consider one experiment; thus, $N_e=1$. First, we run a pilot to estimate the constants \eqref{eq:C1} to \eqref{eq:D3} required for the optimal setting of the DLMC and the DLMC2IS estimators using $N=1000$ outer samples and $M_1=M_2=200$ inner samples. The optimal number of samples and the splitting parameter for a certain tolerance $TOL>0$ are given by \eqref{eq:opt.samples}. For the optimal setting of the MC2LA estimator, we run another pilot using $N=1000$ outer samples to estimate the variance. The results are displayed in Figure~\ref{fig:ex1.opt} as a function of the tolerance $TOL$. For the DLMC estimator, the number of required inner samples $M_1^{\ast}$ is smaller than $M_2^{\ast}$ because the constant $C_1$ only encompasses the variance from the nuisance parameter $\phi$, whereas the constant $C_2$ encompasses the variance from $\theta$ and $\phi$. Importance sampling reduces the number of required inner samples to one, even for small tolerances for the DLMC2IS estimator.

\begin{figure}[ht]
\begin{tikzpicture}
\node[inner sep=0pt] (comp) at (0,0)
    {\includegraphics[width=0.6\textwidth]{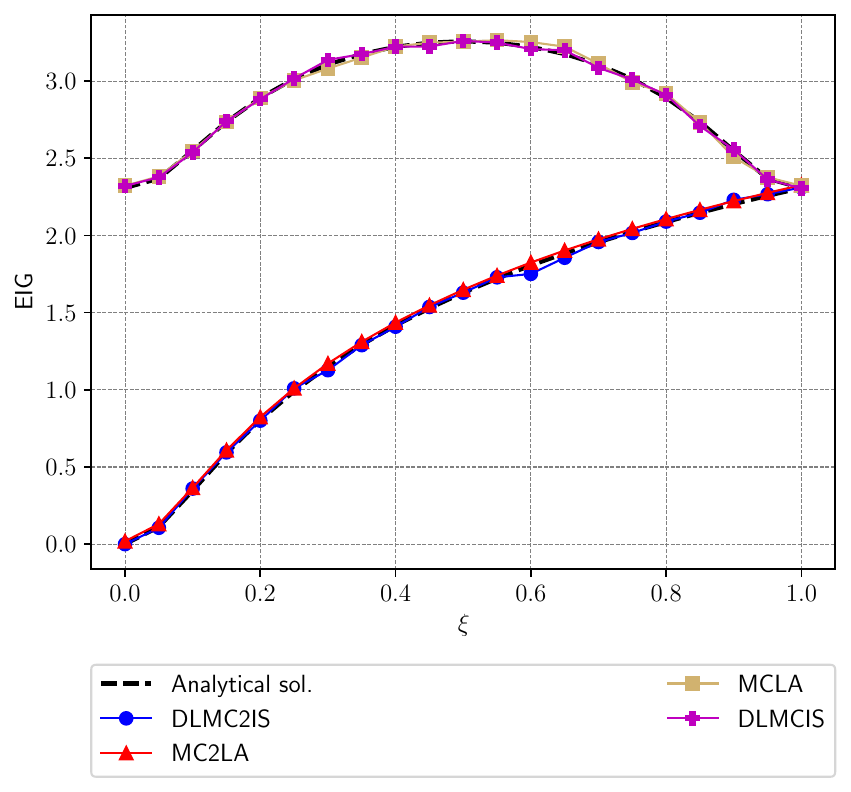}};
		\node at (0.65,2.6) {\small Information gain about $\theta$ and $\phi$};
    \node at (0.65,-1.0) {\small Information gain about $\theta$ only};
\end{tikzpicture}
\caption{Example 1: EIG vs. design parameter $\xi$. Analytical solution (dashed black), DLMC2IS estimator (solid blue), and MC2LA estimator (solid red) for the case with nuisance uncertainty and DLMCIS estimator (solid magenta) and MCLA estimator (solid tan) for the case without nuisance uncertainty.}
\label{fig:ex1.comp}
\end{figure}

We run the estimators for various designs $\xi$ between 0 and 1. The results are presented in Figure~\ref{fig:ex1.comp}. The DLMC2IS and MC2LA estimators indicate that $\xi=1$ is optimal. The model \eqref{eq:ex1.model} is linear; thus, the EIG can be computed analytically. For comparison, we also present the optimal design with no nuisance uncertainty (i.e., both $\theta$ and $\phi$ are considered parameters of interest). For this scenario, we estimate the EIG using the DLMCIS and MCLA estimators developed in \cite{Bec18}, which only use one Laplace approximation. The optimal design is found at $\xi=1/2$. Although a higher overall information gain occurs in this scenario, the information we gain solely about $\theta$ is less at $\xi=1/2$ than at $\xi>1/2$.

\begin{figure}[ht]
	\subfloat[DLMC2IS\label{subfig-1:evt.dlis}]{%
	\begin{tikzpicture}
	\node[inner sep=0pt] (comp) at (0,0)
		{\includegraphics[width=0.45\textwidth]{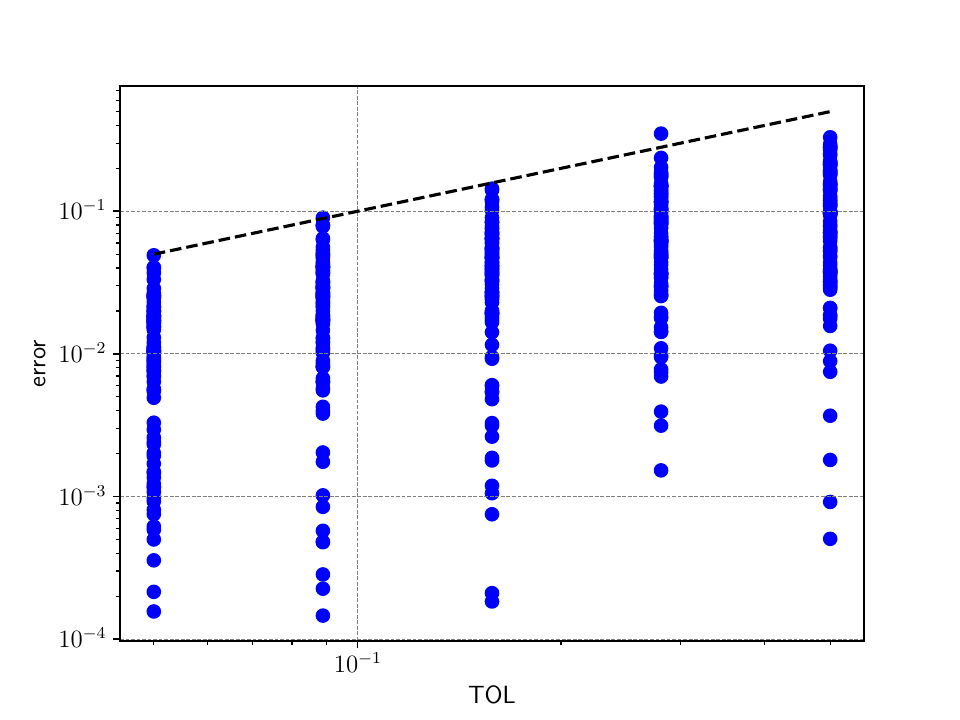}};
		\node at (-2.55,1.1) {\small \red{0}};
    \node at (-1.24,1.4) {\small \red{1}};
     \node at (0.1,1.7) {\small \red{0}};
      \node at (1.4,1.95) {\small \red{1}};
       \node at (2.7,1.95) {\small \red{0}};
		\end{tikzpicture}
	}
	\hfill
	\subfloat[MC2LA\label{subfig-2:evt.mcla}]{%
	\begin{tikzpicture}
	\node[inner sep=0pt] (comp) at (0,0)
		{\includegraphics[width=0.45\textwidth]{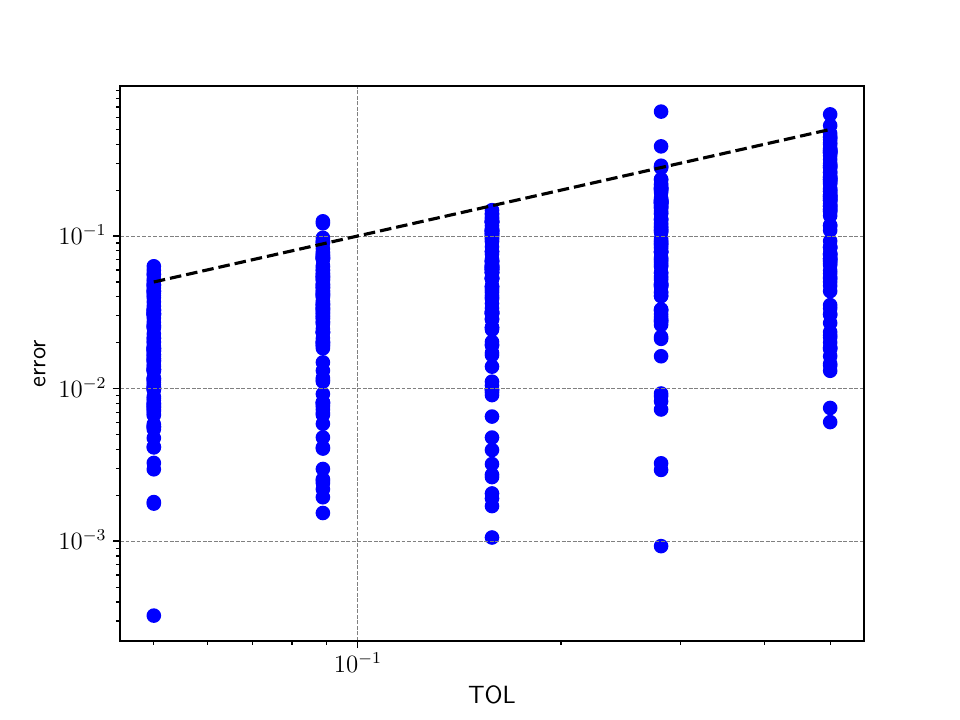}};
		\node at (-2.55,1.1) {\small \red{6}};
    \node at (-1.24,1.4) {\small \red{6}};
     \node at (0.1,1.7) {\small \red{0}};
      \node at (1.2,1.95) {\small \red{4}};
       \node at (2.5,1.95) {\small \red{2}};
		\end{tikzpicture}
	}
	\caption{Example 1: Error vs.~tolerance consistency plot for various tolerances $TOL$ with a predefined confidence parameter $C_\alpha=1.96$ ($\alpha=0.05$). Panel (A) DLMC2IS estimator. Panel (B) MC2LA estimator.}
	\label{fig:ex1.evt}
\end{figure}
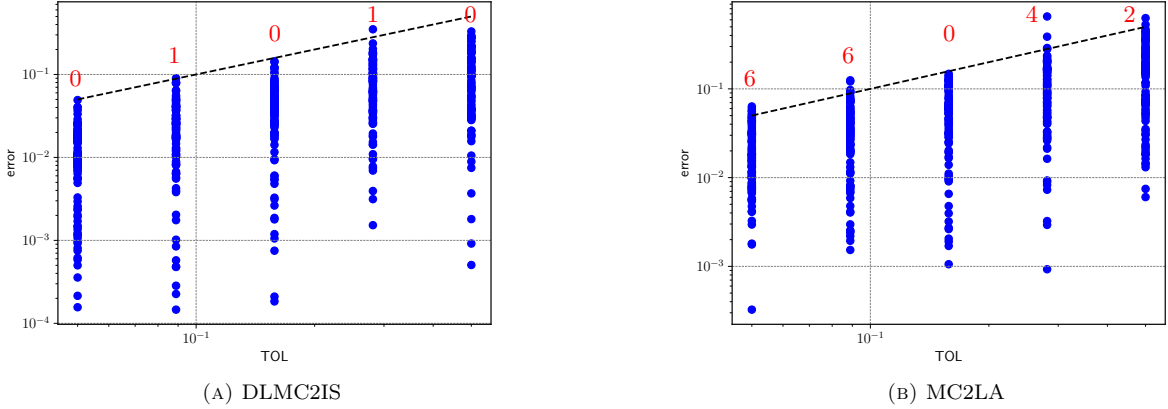

Figure~\ref{fig:ex1.evt} presents 100 runs of the DLMC2IS and MC2LA estimator for various tolerances and design $\xi=1/2$. For every tolerance, the probabilistic error bounds specified by the central limit theorem with confidence constant $C_{\alpha}=1.96$ ($\alpha=0.05$) predict five runs resulting in an error greater than that tolerance. The probabilistic error bound for the DLMC2IS estimator was overly conservative, whereas the probabilistic error bound for the MC2LA estimator was overly optimistic for small tolerances. The reason for the unexpectedly small error of the DLMC2IS estimator is likely due to the fact that even using just one inner sample results in an inner variance that is much smaller than required, thanks to Laplace-based importance sampling.

\subsection{Pharmacokinetics Example}\label{ex:PK}
\begin{figure}[ht]
		\includegraphics[width=0.6\textwidth]{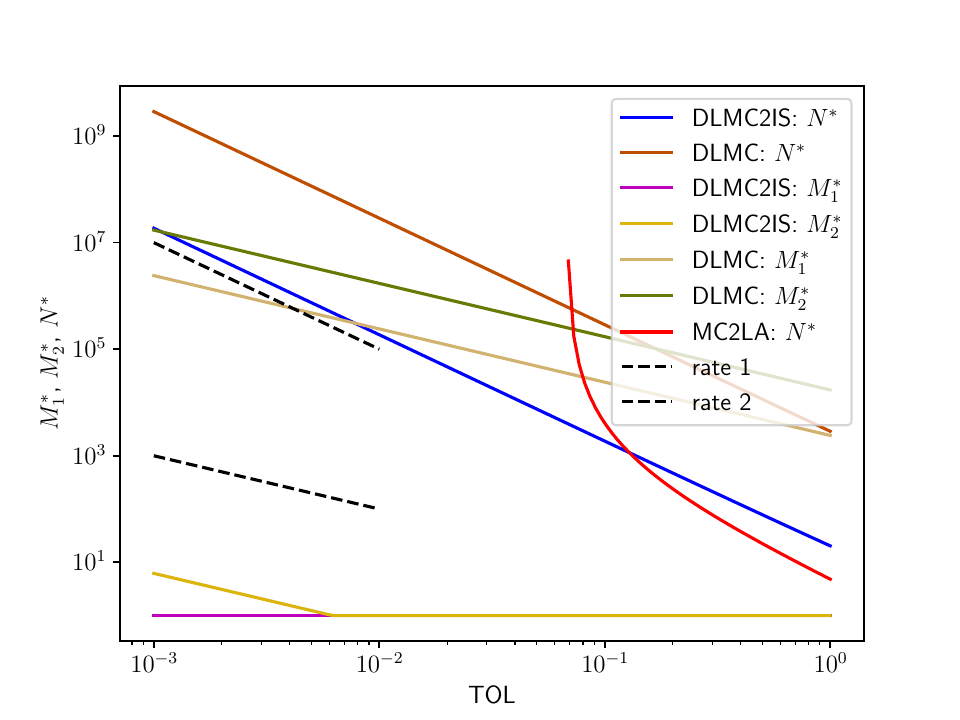}
	\caption{Example 2: Optimal number of outer ($N^{\ast}$) and inner ($M_1^{\ast}$, $M_2^{\ast}$) samples vs. tolerance $TOL$ for the DLMC, DLMC2IS, and MC2LA estimators for geometrically spaced sampling times.}
	\label{fig:ex.pk.opt}
\end{figure}

\begin{figure}[ht]
		\includegraphics[width=0.6\textwidth]{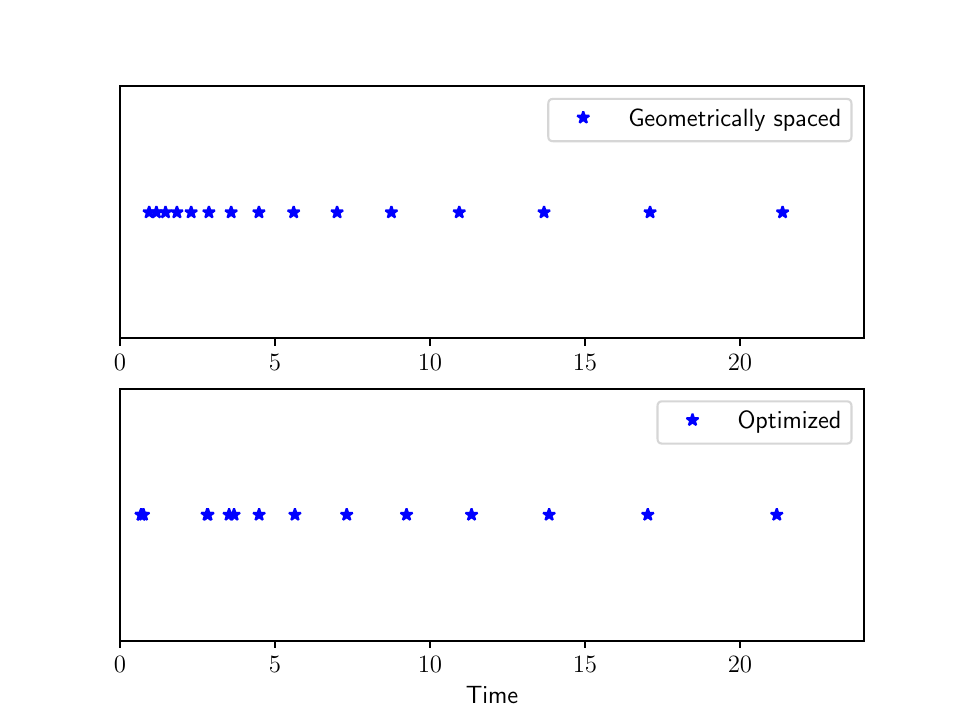}
	\caption{Example 2: Geometrically spaced design \cite{God20} yielding expected information gain (EIG) of 6.12 vs.~optimized design yielding EIG of 6.25. Clustering of measurement times appeared to have a positive impact on the EIG of the experiment. Later measurement times had little effect on the EIG and were mostly unaffected by the optimization.}
	\label{fig:ex.pk.design}
\end{figure}

\begin{figure}[ht]
		\includegraphics[width=0.6\textwidth]{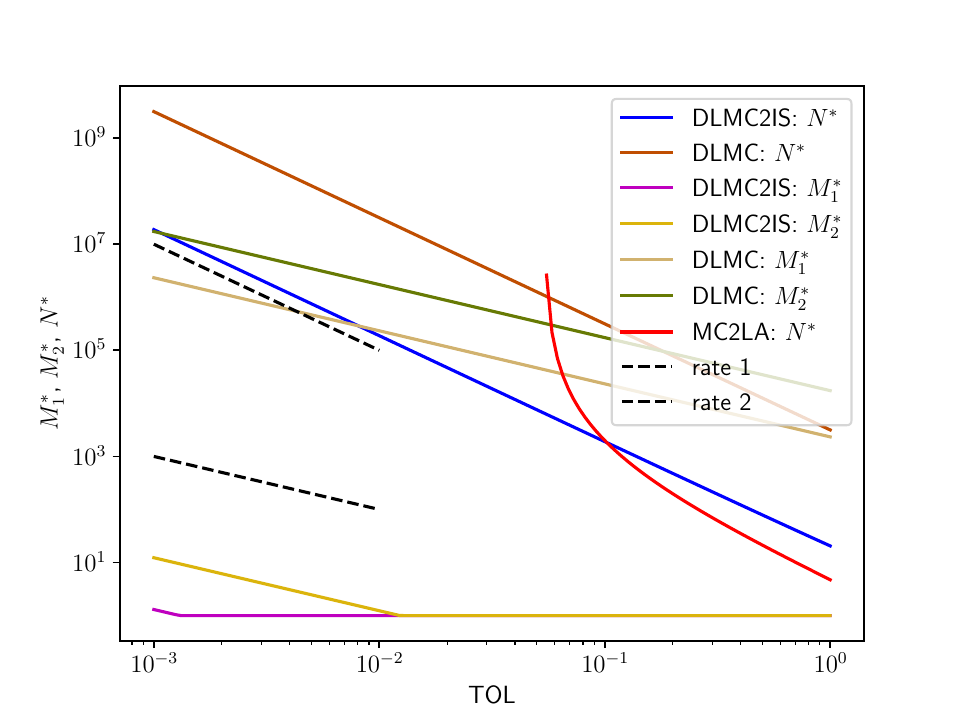}
	\caption{Example 2: Optimal number of outer ($N^{\ast}$) and inner ($M_1^{\ast}$, $M_2^{\ast}$) samples vs. tolerance $TOL$ for the DLMC, DLMC2IS, and MC2LA estimators for optimized sampling times.}
	\label{fig:ex.pk.opt.opt}
\end{figure}

The aim of pharmacokinetics is to learn patient-specific parameters in a medical setting. The following example is based on the setting introduced in \cite{Rya14} and modified in \cite{God20}. After a drug is administered to a patient at time $t_0=0$, a total of 15 blood samples are taken over the next 24 h to determine how fast the drug is absorbed and subsequently eliminated. We follow the simplified approach in \cite{God20}, where only additive noise is present. The design space is $(0,24]^{d_{\xi}}$, where $d_{\xi}=15$. The data model is as follows:
\begin{equation}
    g_j(\bs{\theta},\phi,\bs{\xi})\coloneqq \frac{D}{\phi}\frac{\theta_1}{\theta_1-\theta_2}\left(e^{-\theta_2\xi_j}-e^{-\theta_1\xi_j}\right), \quad 1\leq j \leq 15,
\end{equation}
where $D=400$ indicates the administered dose, $\log(\theta_1)\sim\cl{N}(0,0.05)$ indicates the first-order absorption constant, $\log(\theta_2)\sim\cl{N}(\log(0.1),0.05)$ indicates the first-order elimination constant, and $\log(\phi)\sim\cl{N}(\log(20),0.05)$ indicates the volume of distribution. This last parameter is considered a nuisance parameter in the present work, deviating from the original setting in \cite{God20}. The design $\bs{\xi}=(\xi_1,\ldots,\xi_{15})\coloneqq (t_1,\ldots,t_{15})$ of the experiment signifies the sample times. Moreover, $\bs{\epsilon}\sim\cl{N}(\bs{0},\bs{\Sigma_{\bs{\epsilon}}})$, where
\begin{equation}
    \bs{\Sigma_{\bs{\epsilon}}}=\begin{pmatrix}
        10^{-2}&0&\cdots&0\\
        0&\ddots&\ddots&\vdots\\
        \vdots&\ddots&\ddots&0\\
        0&\cdots&0&10^{-2}
    \end{pmatrix}
\end{equation}
signifies the observation noise. The work \cite{God20} found that a geometrically spaced design $\xi_j=0.94\times1.25^{j-1}$, where $1\leq j\leq 15$, performed best for their experiments, which did not consider nuisance uncertainty. We adopt this choice as a starting point for our optimization and present the optimal number of samples for the DLMC, DLMC2IS, and MC2LA estimators in Figure~\ref{fig:ex.pk.opt}. For the pilot of the DLMC estimator, we used $N=300$ outer samples and $M_1=M_2=13000$ inner samples. For the pilot of the DLMC2IS estimator, we used $N=2000$ outer samples and $M_1=M_2=200$ inner samples. Finally, for the pilot of the MC2LA estimator, we used $N=2000$ outer samples and the DLMC2IS estimator to ascertain the bias resulting from the Laplace approximation.
To improve upon this design, we used a greedy minibatch stochastic gradient descent algorithm, where $N=300$ outer samples of the MC2LA estimator are used to optimize the EIG with respect to $\xi_1$ with $(\xi_2,\ldots,\xi_{15})$ fixed. Next, the EIG is optimized with respect to $\xi_2$ with $(\xi_1,\xi_3,\ldots,\xi_{15})$ fixed, and so on and so forth. That is, we employ a combination of stochastic gradient descent and coordinate descent methods to find a local minimum of the EIG. The geometrically spaced design is used as a starting point. A comparison between both design choices is presented in Figure~\ref{fig:ex.pk.design}. Optimization of measurement times indicates that a certain clustering is beneficial for the EIG, whereas measurements towards the end of the 24 h appeared to have little effect on the EIG. A clustering effect of sampling times was also observed in \cite{Rya14}.
The optimal number of inner and outer samples for the optimized design is displayed in Figure~\ref{fig:ex.pk.opt.opt}. The design had minimal effect on the optimal number of samples; however, the EIG was improved for the optimized design, as the DLMC2IS estimator with a tolerance of $TOL=10^{-2}$ yielded an $EIG$ of 6.25 for the optimized design and an EIG of 6.12 for the geometrically spaced design.

\subsection{Electrical Impedance Tomography Example I}\label{ex:EIT.1}
For this example, we consider a more challenging model based on the solution operator of a PDE. Solving the PDE for EIT is generally not possible in closed form; therefore, we employ a finite element method (FEM) approximation instead. This example demonstrates the practical applicability of the derived estimators and their ability to incorporate approximate models $\bs{g}_h$, where $h$ is the mesh-discretization parameter. We examine a two-dimensional model of a composite laminate material with a fiber structure. The material consists of two plies, both conducting an electric current more easily along the direction of their fibers than transversal to them. The experimental setup involves attaching five electrodes to the top and five to the bottom boundaries, injecting, and measuring electric current. The experimenter can learn about the fiber angles in each ply; hence, these angles are considered parameters of interest. The electrode positions can be chosen freely in principle but affect the meaningfulness of the measurements. Hence, this is the design the experimenter aims to optimize. We further assume that the exact conductivities of the plies are known only up to a concentrated normal distribution. Therefore, they are considered nuisance parameters. We consider a rectangular domain $\cl{D}=\cl{D}_1\cup\cl{D}_2=[0,20]\times[0,1]\cup [0,20]\times[1,2]$, consisting of two subdomains. For the quasi-static potential field $u$, current flux $\bs{j}$, and conductivity field $\bs{\bar{\sigma}}$, we solve the PDE
\begin{align}
    \nabla\cdot\bs{j}(\bs{x},\omega){}&=0\quad \text{in $\cl{D}$ and}\\
    \bs{j}(\bs{x},\omega){}&=\bs{\bar{\sigma}}(\omega)\cdot\nabla u(\bs{x},\omega),
\end{align}
where $\bs{x}\in\cl{D}$,
\begin{equation}
    \bs{\bar{\sigma}}(\omega)=\bs{Q}(\theta_i(\omega))^{\trans}\cdot\bs{\sigma}(\phi_i(\omega))\cdot\bs{Q}(\theta_i(\omega)), \quad i=1,2,
\end{equation}
\begin{equation}
    \bs{Q}(\theta_i)=\begin{pmatrix}
        \cos(\theta_i) & 0 & -\sin(\theta_i)\\ 0 & 1 & 0\\ \sin(\theta_i) & 0 & \cos(\theta_i)
    \end{pmatrix},\quad i=1,2,
\end{equation}
and
\begin{equation}
    \bs{\sigma}(\phi_i)=\begin{pmatrix}
        \sigma_1(\phi_i) & 0 & 0\\ 0 & \sigma_2(\phi_i) & 0\\ 0 & 0&\sigma_3(\phi_i)
    \end{pmatrix}, \quad i=1,2.
\end{equation}
For the fiber angle in $\cl{D}_1$, we assume the following prior distribution:
\begin{equation}
    \theta_1\sim\pi(\theta_1)=\cl{U}\left(-\frac{\pi}{4}-0.05,-\frac{\pi}{4}+0.05\right)
\end{equation}
and for the fiber angle in $\cl{D}_2$, we assume the prior distribution:
\begin{equation}
    \theta_2\sim\pi(\theta_2)=\cl{U}\left(\frac{\pi}{4}-0.05,\frac{\pi}{4}+0.05\right).
\end{equation}
For the conductivities, we consider the prior distribution:
\begin{equation}
 \sigma_j(\phi_i)=\exp(\mu_j+\phi_i),\quad j=1,2,3,\quad i=1,2,   
\end{equation}
where
\begin{equation}
    \phi_i=\sigma_{\phi}z_i, \quad z_i\stackrel{\mathrm{iid}}{\sim}\cl{N}(0,1), \quad i=1,2.
\end{equation}
Moreover, $\mu_1=\log(0.1)$ and $\mu_2=\mu_3=\log(0.02)$ and $\sigma^2_{\phi}$ is the covariance of $\phi_i$ for all $i=1,2$. The design $\bs{\xi}=(\xi_1,\xi_2)=[0,2]\times[0,2]$ signifies a shift between top and bottom electrodes and the distance between electrodes, respectively. For a detailed description of the problem setting, including boundary conditions, and a finite element formulation, see \cite{Bar22, Bec18, Som92}.

\begin{figure}[ht]
		\includegraphics[width=0.6\textwidth]{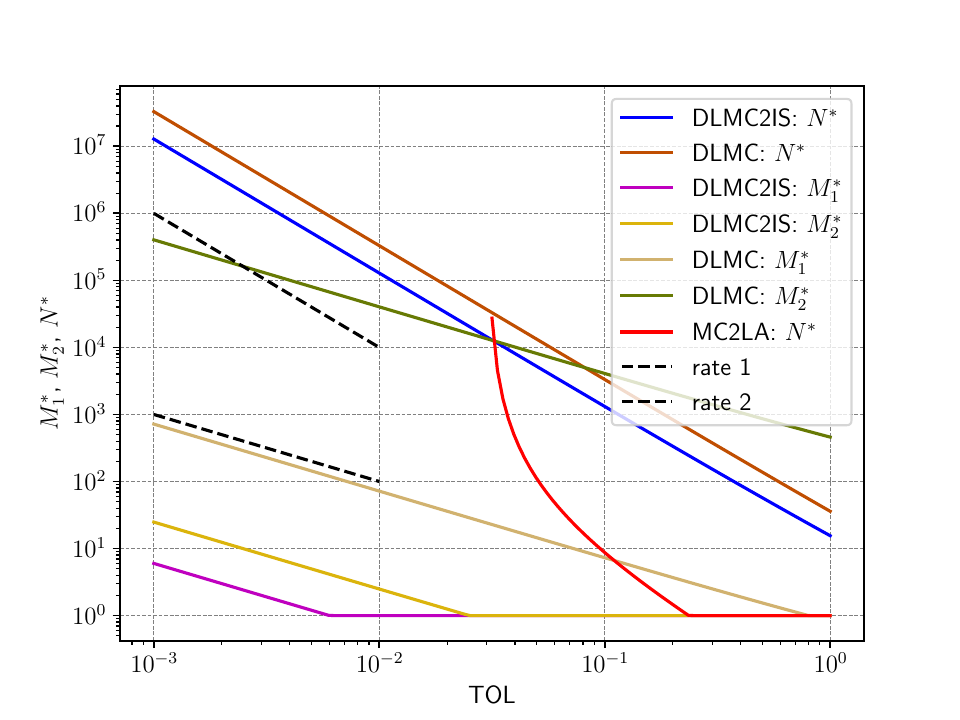}
	\caption{Example 3: Optimal number of outer ($N^{\ast}$) and inner ($M_1^{\ast}$, $M_2^{\ast}$) samples vs. tolerance $TOL$ for the DLMC, DLMC2IS, and MC2LA estimators.}
	\label{fig:ex2.opt}
\end{figure}

\begin{figure}[ht]
		\includegraphics[width=0.6\textwidth]{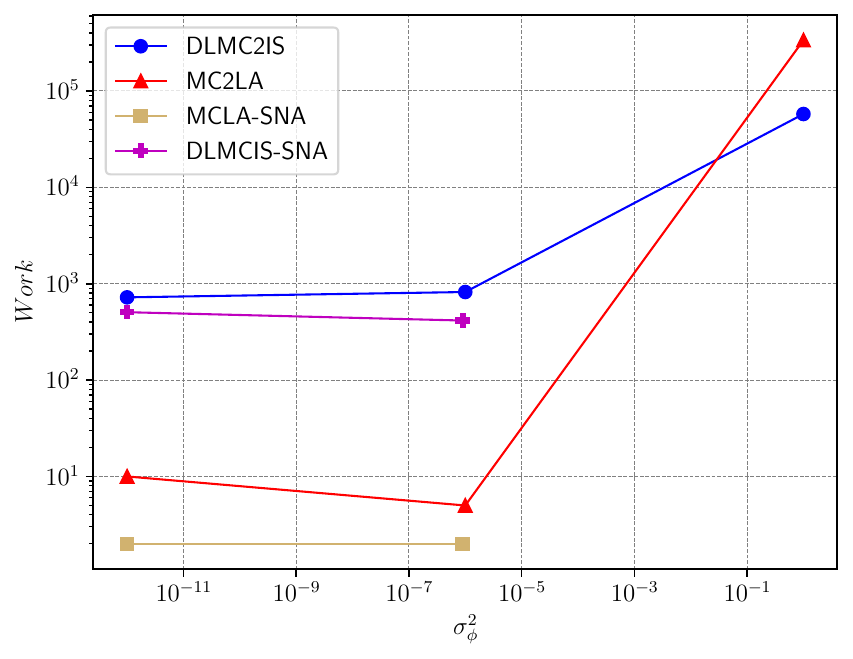}
	\caption{Example 3: Computational work vs.~nuisance covariance $\sigma_{\phi}^2$ for the MCLA with small-noise approximation, DLMCIS with small-noise approximation, DLMC2IS, and MC2LA estimators.}
	\label{fig:ex2.work}
\end{figure}

As in the previous example, we start by running a pilot with $N=50$ outer samples for the MC2LA estimator and $N=50$ outer samples, with $M_1=M_2=10$ inner samples for the DLMC2IS and DLMC estimators. The bias resulting from the Laplace approximations in the MC2LA estimator is measured by comparison with the DLMC2IS estimator results. In addition, we estimate the FEM constants $\eta$, $C_{\rm{disc}}$, and $\gamma$. The optimal number of samples for these estimators is depicted in Figure~\ref{fig:ex2.opt} as a function of the error tolerance $TOL$. Next, we compare the computational work required for the MC2LA and DLMC2IS estimators and for the MCLA and DLMCIS estimators with small-noise approximation developed in \cite{Bar22} as a function of the covariance of the nuisance parameters in Figure~\ref{fig:ex2.work}. The computational work for the MC2LA and the MCLA with small-noise approximation estimators is given by $N\times h^{-\gamma}$. The computational work for the DLMC2IS estimator is given by $N\times(M_1+M_2)\times h^{-\gamma}$. Finally, the computational work for the DLMCIS estimator with small-noise approximation is given by $N\times M\times h^{-\gamma}$. The small-noise approximation is only applicable up to a small nuisance covariance $\sigma_{\phi}^2\approx 9.03\times 10^{-7}$. However, the cost for the MC2LA and the DLMC2IS estimators increase sharply for a larger nuisance covariance.

\subsection{Electrical Impedance Tomography Example II}
\begin{figure}[ht]
	\subfloat[Current streamlines for design $\bs{\xi}=(0.5,2.0)$.\label{subfig-1:currents.1}]{%
		\includegraphics[width=0.95\textwidth,trim={0cm 20cm 0cm 20cm},clip]{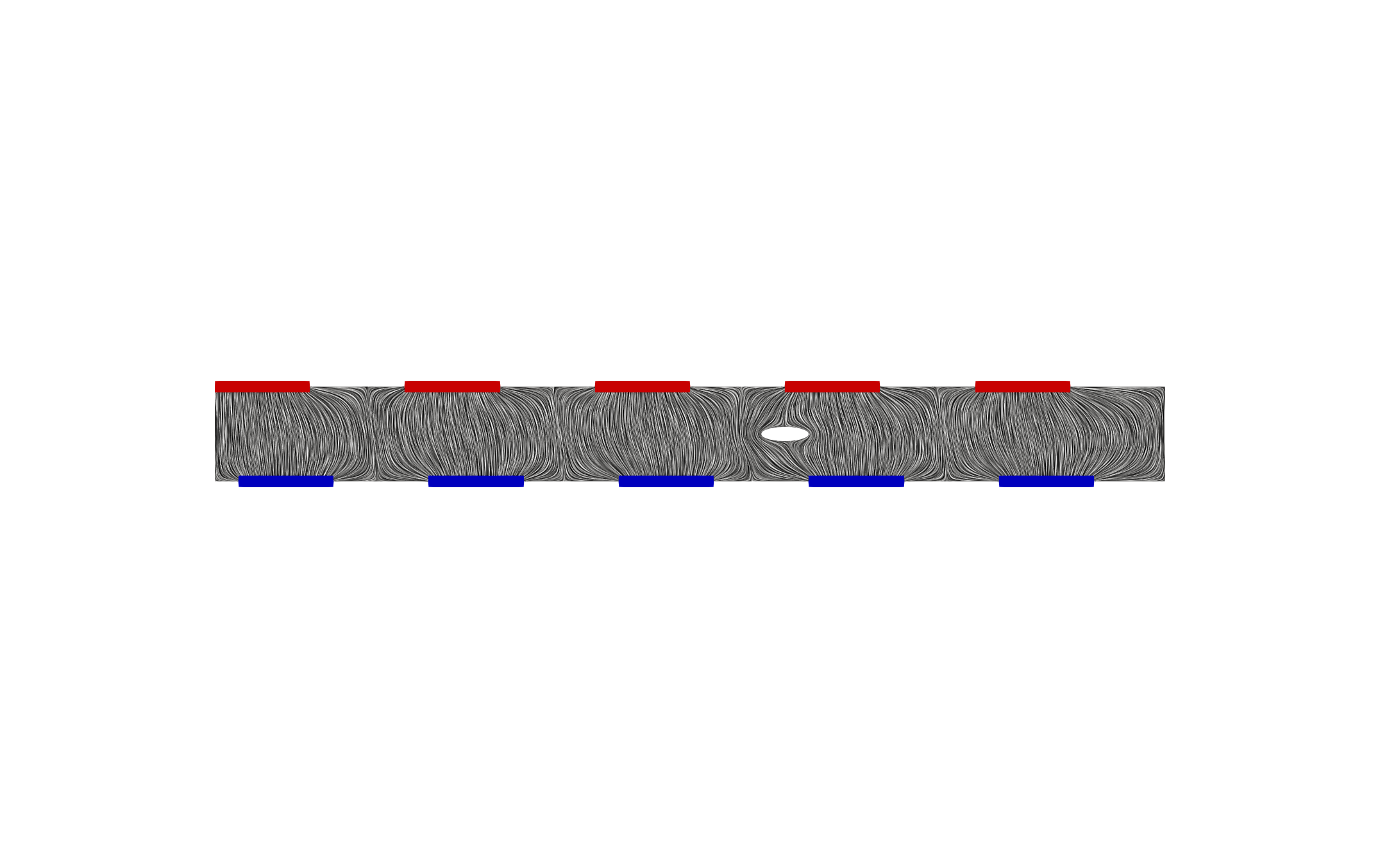}
	}
	\\
	\subfloat[Current streamlines for design $\bs{\xi}=(2.0,0.5)$.\label{subfig-2:currents.2}]{%
		\includegraphics[width=0.95\textwidth,trim={0cm 20cm 0cm 20cm},clip]{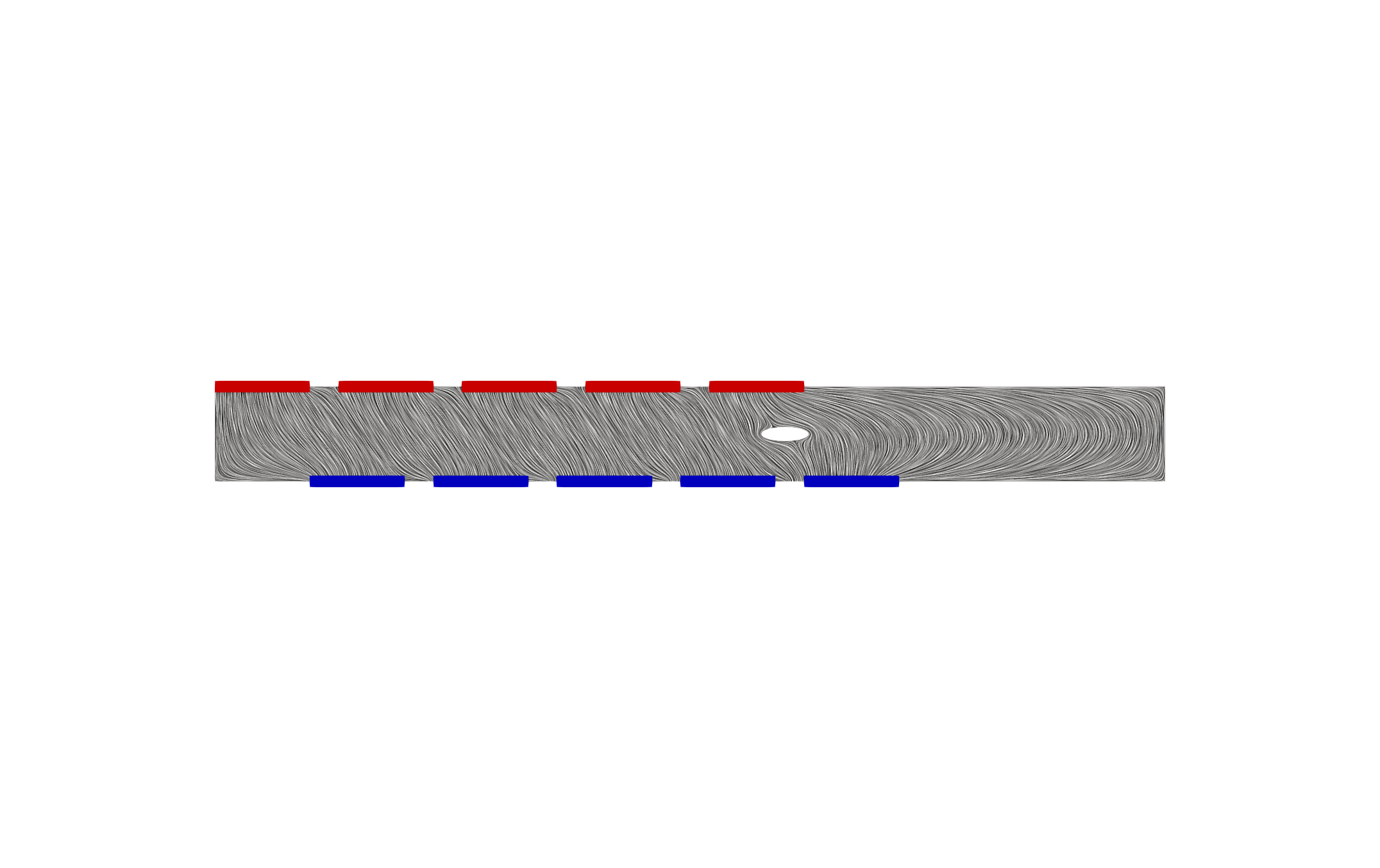}
	}
	\caption{Example 4: Inlet: red electrodes. Outlet: blue electrodes. The current flow is affected by the electrode positions.}
	\label{fig:ex.3.xi}
\end{figure}

\begin{figure}[ht]
	\subfloat[Expected information gain regarding the fiber angles and the vertical axis of the ellipsoid exclusion.\label{subfig-1:xi.3}]{%
		\includegraphics[width=0.45\textwidth,trim={0cm 0cm 4cm 0cm},clip]{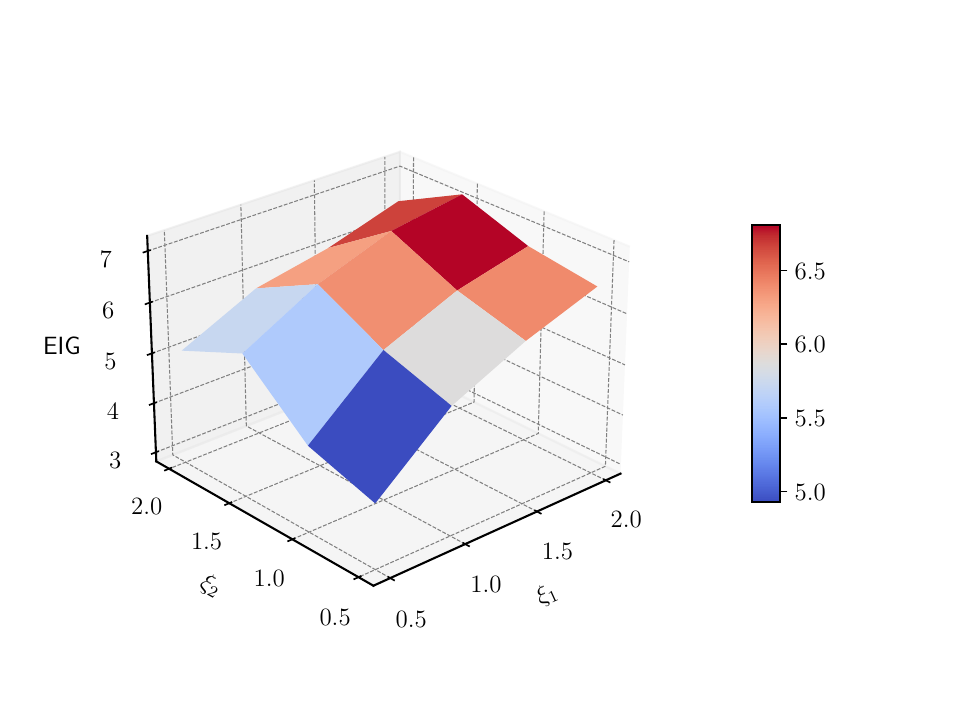}
	}
	\subfloat[Expected information gain regarding the fiber angles only.\label{subfig-2:xi.4}]{%
		\includegraphics[width=0.45\textwidth,trim={0cm 0cm 4cm 0cm},clip]{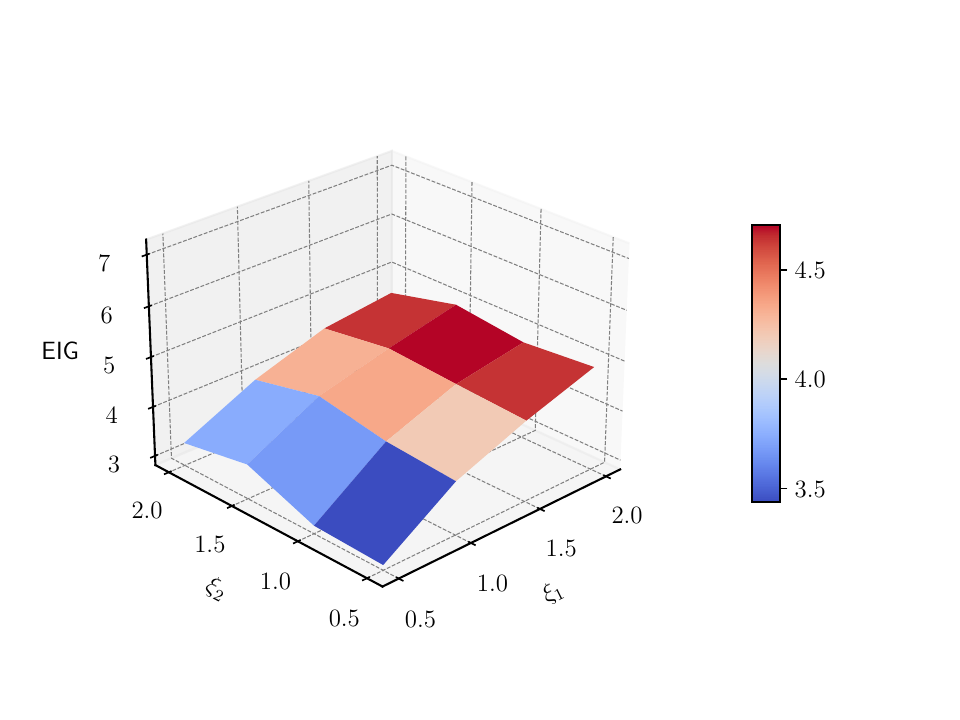}
	}
	\caption{Example 4: Expected information gain (EIG) as a function of the design $\bs{\xi}$. The response surface changes when considering $\bs{\theta}$ and $\phi$ as parameters of interest (Panel~\protect\subref{subfig-1:xi.3}), or $\bs{\theta}$ only (Panel~\protect\subref{subfig-2:xi.4}). The EIG regarding the fiber angles is virtually unaffected by the distance between the electrodes ($\xi_2$) in this case.}
	\label{fig:ex.3.xi.3}
\end{figure}

Next, we consider a slightly different setup from the previous example. The electric conductivity is now fixed at $\sigma_1=0.1$, $\sigma_2=\sigma_3= 0.02$, and no longer a nuisance parameter. Furthermore, we assume that a small ellipsoid exclusion exists between the two plies. The vertical axis of this ellipsoid is the new nuisance parameter, with the following distribution:
\begin{equation}
  \phi\sim\cl{U}(0.2, 0.4).
\end{equation}
The horizontal axis is considered fixed at 1 and the center is fixed at (12,1). The parameters of interest are still the angles of the fibers in each ply. The FEM formulation remains the same as in the previous example, except for the different mesh that now incorporates the ellipsoid hole (Figure~\ref{fig:ex.3.xi}).
Figure~\ref{fig:ex.3.xi.3} demonstrates that although the design $\bs{\xi}=(2.0,1.5)$ is optimal regardless of whether $\bs{\theta}$ and $\phi$ or only $\bs{\theta}$ are considered parameters of interest, the overall response surface changes depending on the parameters that are considered to be of interest. In particular, $\xi_2$, the distance between electrodes, is almost irrelevant when recovering the fiber angles, whereas it has a much more prominent effect when recovering the height of the exclusion as well. The design choices $\bs{\xi}=(2.0,1.0)$ and $\bs{\xi}=(2.0,0.5)$ yield EIG that is within the tolerance $TOL=0.2$ of the optimal choice when considering nuisance parameters. The MCLA estimator with optimal sampling was used to estimate the EIG without nuisance uncertainty. The pilot run for this estimator was performed using $N=150$ samples. The DLMC2IS estimator was used with optimal sampling to estimate the EIG with nuisance uncertainty. For the pilot run, we used $N=M_1=M_2=200$ samples. The bias introduced by the Laplace approximation \eqref{ptilde} rendered the MC2LA estimator ineffective for this example for a tolerance of $TOL=0.2$; thus the DLMC2IS estimator was used instead.

\section{Conclusion}
We propose two estimators for the expected information gain under nuisance uncertainty: the Monte Carlo double-Laplace estimator based on Laplace's method for the nuisance parameters and the Laplace approximation for the parameters of interest, and the double-loop Monte Carlo double importance sampling estimator based on two Laplace approximations. We demonstrate the applicability of these estimators in four numerical examples, showcasing their computational efficiency provided by Laplace-based integral approximations.

\section{Statements and Declarations}
\subsection{Acknowledgments}
This publication is based upon work supported by the King Abdullah University of Science and Technology (KAUST) Office of Sponsored Research (OSR) under Award No.~OSR-2019-CRG8-4033, the Alexander von Humboldt Foundation, the Deutsche Forschungsgemeinschaft (DFG, German Research Foundation) -- 333849990/GRK2379 (IRTG Hierarchical and Hybrid Approaches in Modern Inverse Problems), and was partially supported by the Flexible Interdisciplinary Research Collaboration Fund at the University of Nottingham Project ID 7466664.

\subsection{Competing interests}
Authors Luis Espath and Ra\'{u}l Tempone are associate editors of the Statistics and Computing Journal.

\subsection{Data availability}
The data that support the findings of this study are available from the corresponding author upon reasonable request.

%-------------------------------------------------------------------------------%

\appendix

\section{Error estimates for the double-loop Monte Carlo estimator with two inner loops}\label{ap:bias}
\subsection*{Derivation of the bias error approximation}
Following \cite{Bec18}, we derive an estimate for the bias, given as follows:
\begin{equation}\label{eq:bias}
|\mathbb{E}[I_{\rm{DL}}]-I|.
\end{equation}
We recall that $I_{\rm{DL}}$ is given as follows:
\begin{align}\label{eq:Dkl.sampled}
I_{\rm{DL}}={}& \frac{1}{N}\sum_{n=1}^N\log\left(\frac{1}{M_1}\sum_{m=1}^{M_1}p(\bs Y^{(n)}|\bs{\theta}^{(n)},\bs{\varphi}^{(n,m)})\right)-\log\left(\frac{1}{M_2}\sum_{k=1}^{M_2}p(\bs Y^{(n)}|\bs{\vartheta}^{(n,k)},\bs{\varphi}^{(n,k)})\right),\nonumber\\
\coloneqq{}&\frac{1}{N}\sum_{n=1}^N\log(\hat{p}_{M_1}(\bs{Y}^{(n)}|\bs{\theta}^{(n)}))-\log(\hat{p}_{M_2}(\bs{Y}^{(n)})).
\end{align}
Replacing \eqref{eq:Dkl.sampled} in \eqref{eq:bias} yields
\begin{align}
|\mathbb{E}[I_{\rm{DL}}]-I|={}&\left|\mathbb{E}\left[\frac{1}{N}\sum_{n=1}^N\log\left(\frac{\hat{p}_{M_1}(\bs Y^{(n)}|\bs{\theta}^{(n)})} {\hat{p}_{M_2}( \bs Y^{(n)}) } \right)\right]-\mathbb{E}\left[\log\left(\frac{p(\bs Y|\bs{\theta})} {p( \bs Y) } \right)\right]\right|,\nonumber\\[4pt]
={}&|\underbrace{\mathbb{E}[\log(\hat{p}_{M_1}(\bs Y|\bs{\theta}))]}_{I}-\underbrace{\mathbb{E}[\log(p(\bs Y|\bs{\theta}))]}_{II}-\underbrace{\mathbb{E}[\log(\hat{p}_{M_2}( \bs Y))]}_{III}+\underbrace{\mathbb{E}[\log(p(\bs Y))]}_{IV}|.
\end{align}

The estimation centers on the following Taylor expansion of $\log(X)$ around $\mathbb{E}[X]$:
\begin{equation}\label{tayl}
\log(X) = \log(\mathbb{E}[X]) + \frac{1}{\mathbb{E}[X]}(X-\mathbb{E}[X])-\frac{1}{2}\frac{1}{\mathbb{E}[X]^2}(X-\mathbb{E}[X])^2+\int_0^1\frac{1-s}{\left(\bb{E}[X]+s\left(X-\bb{E}[X]\right)\right)^3}\di{}s(X-\mathbb{E}[X])^3.
\end{equation}

As $\hat{p}_{M_2}(\bs Y)$ is an unbiased estimator, we have
\[\mathbb{E}[\hat{p}_{M_2}(\bs Y)|\bs{Y}]=p(\bs Y).\]
By the Taylor expansion and taking the expected value, we obtain
\begin{align}\label{eq:Taylor.M2}
\mathbb{E}[\log(\hat{p}_{M_2}(\bs Y))|\bs{Y}]=\mathbb{E}[\log(p(\bs Y))|\bs{Y}]{}&-\frac{1}{2}\mathbb{E}\left[\frac{1}{p^2(\bs Y)}(\hat{p}_{M_2}(\bs Y)-p(\bs Y))^2|\bs{Y}\right]\nonumber\\
{}&+\bb{E}\left[\int_0^1\frac{1-s}{\left(p(\bs Y)+s\left(\hat{p}_{M_2}(\bs Y)-p(\bs Y)\right)\right)^3}\di{}s(\hat{p}_{M_2}(\bs Y)-p(\bs Y))^3\Bigg|\bs{Y}\right],
\end{align}
where the linear term vanishes because of the expected value. The quadratic term can be rewritten as follows:
\begin{align}
\mathbb{E}\left[\frac{1}{p^2(\bs Y)}(\hat{p}_{M_2}(\bs Y)-p(\bs Y))^2|\bs{Y}\right]={}&\frac{\mathbb{V}[\hat{p}_{M_2}(\bs Y)|\bs Y]}{p^2(\bs Y)},\\
={}&\frac{\mathbb{V}\left[\frac{1}{M_2}\sum_{k=1}^{M_2}p(\bs Y|\bs{\vartheta}^{(k)},\bs{\varphi}^{(k)})|\bs Y\right]}{p^2(\bs Y)},\nonumber\\
={}&\frac{1}{M_2}\frac{\mathbb{V}\left[p(\bs Y|\bs{\vartheta},\bs{\varphi})\Big|\bs Y\right]}{p^2(\bs Y)}.
\end{align}
The third-order term has the following bound:
\begin{align}\label{eq:holder}
    \bb{E}{}&\left[\int_0^1\frac{1-s}{\left(p(\bs Y)+s\left(\hat{p}_{M_2}(\bs Y)-p(\bs Y)\right)\right)^3}\di{}s(\hat{p}_{M_2}(\bs Y)-p(\bs Y))^3\Bigg|\bs{Y}\right]\nonumber\\
    {}&\leq\bb{E}\left[\left(\int_0^1\frac{1-s}{\left(p(\bs Y)+s\left(\hat{p}_{M_2}(\bs Y)-p(\bs Y)\right)\right)^3}\di{}s\right)^{p}\Bigg|\bs{Y}\right]^{\frac{1}{p}}\bb{E}\left[(\hat{p}_{M_2}(\bs Y)-p(\bs Y))^{3q}\Bigg|\bs{Y}\right]^{\frac{1}{q}}
\end{align}
by H\"{o}lder's inequality for some $1\leq p,q\leq\infty$, where $1/p+1/q=1$. For the last term in \eqref{eq:holder}, we have by the discrete Burkholder--Davis--Gundy inequality \cite{Bur72, Gil19} that:
\begin{align}
    \bb{E}\left[(\hat{p}_{M_2}(\bs Y)-p(\bs Y))^{3q}\Bigg|\bs{Y}\right]^{\frac{1}{q}}{}&=\bb{E}\left[\left(\frac{1}{M_2}\sum_{k=1}^{M_2}p(\bs Y|\bs{\vartheta}^{(k)},\bs{\varphi}^{(k)})-p(\bs Y)\right)^{3q}\Bigg|\bs{Y}\right]^{\frac{1}{q}},\nonumber\\
    {}&\leq \frac{C_{\rm{BDG}}}{M_2^{\frac{3}{2}}}\bb{E}\left[\left(p(\bs Y|\bs{\vartheta},\bs{\varphi})-p(\bs Y)\right)^{3q}\Bigg|\bs{Y}\right]^{\frac{1}{q}}
\end{align}
for some constant $C_{\rm{BDG}}$ for almost all $\bs{Y}$. From \eqref{eq:Taylor.M2}, we obtain that
\begin{equation}\label{ed}
  \underbrace{\mathbb{E}[\log(\hat{p}_{M_2}(\bs Y))]}_{III}=\underbrace{\mathbb{E}[\log(p(\bs Y))]}_{IV}-\frac{1}{2}\frac{1}{M_2}\mathbb{E}\left[\frac{\mathbb{V}\left[p(\bs Y|\bs{\vartheta},\bs{\varphi})\Big|\bs Y\right]}{p^2(\bs Y)}\right]+\mathcal{O}_{\bb{P}}\left(\frac{1}{M_2^{\frac{3}{2}}}\right).
\end{equation}
Similarly, we obtain
\begin{equation}
\mathbb{E}[\hat{p}_{M_1}(\bs Y|\bs{\theta})|\bs{Y},\bs{\theta}]=p(\bs Y|\bs{\theta}),
\end{equation}
and
\begin{equation}
\mathbb{E}[\log(\hat{p}_{M_1}(\bs Y|\bs{\theta}))|\bs{Y},\bs{\theta}]=\mathbb{E}[\log(p(\bs Y|\bs{\theta}))|\bs{Y},\bs{\theta}]-\frac{1}{2}\mathbb{E}\left[\frac{1}{p^2(\bs Y|\bs{\theta})}(\hat{p}_{M_1}(\bs Y|\bs{\theta})-p(\bs Y|\bs{\theta}))^2|\bs{Y},\bs{\theta}\right]+\mathcal{O}_{\bb{P}}\left(\frac{1}{M_1^{\frac{3}{2}}}\right),
\end{equation}
where the quadratic term can be written as follows:
\begin{align}
\mathbb{E}\left[\frac{1}{p^2(\bs Y|\bs{\theta})}(\hat{p}_{M_1}(\bs Y|\bs{\theta})-p(\bs Y|\bs{\theta}))^2|\bs{Y},\bs{\theta}\right]={}&\frac{\mathbb{V}[\hat{p}_{M_1}(\bs Y|\bs{\theta})|\bs Y,\bs{\theta}]}{p^2(\bs Y|\bs{\theta})},\nonumber\\
={}&\frac{\mathbb{V}\left[\frac{1}{M_1}\sum_{m=1}^{M_1}p(\bs Y|\bs{\theta},\bs{\varphi}^{(m)})|\bs Y,\bs{\theta}\right]}{p^2(\bs Y|\bs{\theta})},\nonumber\\
={}&\frac{1}{M_1}\frac{\mathbb{V}\left[p(\bs Y|\bs{\theta},\bs{\varphi})\Big|\bs Y,\bs{\theta}\right]}{p^2(\bs Y|\bs{\theta})}.
\end{align}
This calculation results in
\begin{equation}\label{ed2}
  \underbrace{\mathbb{E}[\log(\hat{p}_{M_1}(\bs Y|\bs{\theta}))]}_{I}=\underbrace{\mathbb{E}[\log(p(\bs Y|\bs{\theta}))]}_{II}-\frac{1}{2}\frac{1}{M_1}\mathbb{E}\left[\frac{\mathbb{V}\left[p(\bs Y|\bs{\theta},\bs{\varphi})\Big|\bs Y,\bs{\theta}\right]}{p^2(\bs Y|\bs{\theta})}\right]+\mathcal{O}_{\bb{P}}\left(\frac{1}{M_1^{\frac{3}{2}}}\right).
\end{equation}
After combining everything, we obtain
\begin{equation}
  |\mathbb{E}[I_{\rm{DL}}]-I|\approx\left|\frac{1}{2M_2}\mathbb{E}\left[\frac{\mathbb{V}[p(\bs Y|\bs{\vartheta},\bs{\varphi})|\bs Y]}{p^2(\bs Y)}\right]-\frac{1}{2M_1}\mathbb{E}\left[\frac{\mathbb{V}[p(\bs Y|\bs{\theta},\bs{\varphi})|\bs Y,\bs{\theta}]}{p^2(\bs Y|\bs{\theta})}\right]
\right|.
\end{equation}

\subsection*{Derivation of the statistical error approximation}
For the variance estimation, we obtain
\begin{equation}
\mathbb{V}[I_{\rm{DL}}]=\mathbb{V}\left[\frac{1}{N}\sum_{n=1}^N\log(\hat{p}_{M_1}(\bs{Y}^{(n)}|\bs{\theta}^{(n)}))-\log(\hat{p}_{M_2}(\bs{Y}^{(n)}))\right].
\end{equation}
By the law of total variance, we can rewrite this as
\begin{align}\label{ve}
\mathbb{V}[I_{\rm{DL}}]={}&\frac{1}{N}\mathbb{V}\left[\mathbb{E}\left[\log(\hat{p}_{M_1}(\bs Y|\bs{\theta}))-\log(\hat{p}_{M_2}(\bs Y))|\bs Y, \bs \theta\right]\right]\\\label{ev}
{}&+\frac{1}{N}\mathbb{E}\left[\mathbb{V}\left[\log(\hat{p}_{M_1}(\bs Y|\bs{\theta}))-\log(\hat{p}_{M_2}(\bs Y))|\bs Y, \bs \theta\right]\right].
\end{align}

Using \eqref{ed} and \eqref{ed2}, we can rewrite \eqref{ve} as
\begin{multline}
\frac{1}{N}\mathbb{V}\left[\mathbb{E}\left[\log(\hat{p}_{M_1}(\bs Y|\bs{\theta}))|\bs Y,\bs{\theta}\right]-\mathbb{E}\left[\log(\hat{p}_{M_2}(\bs Y))|\bs Y\right]\right]\nonumber\\[4pt]
=\frac{1}{N}\mathbb{V}\left[\log\left(\frac{p(\bs Y|\bs{\theta})}{p(\bs Y)}\right)+\frac{1}{2M_2}\frac{\mathbb{V}[p(\bs Y|\bs{\vartheta},\bs{\varphi})|\bs Y]}{p^2(\bs Y)}
-\frac{1}{2M_1}\frac{\mathbb{V}[p(\bs Y|\bs{\theta},\bs{\varphi})|\bs Y,\bs{\theta}]}{p^2(\bs Y|\bs{\theta})}\right]\nonumber\\[4pt]
+\mathcal{O}_{\bb{P}}\left(\frac{1}{NM_1^2}\right)+\mathcal{O}_{\bb{P}}\left(\frac{1}{NM_2^2}\right),
\end{multline}
which results in
\begin{align}\label{eq:cov.terms}\nonumber
\frac{1}{N}{}&\mathbb{V}\left[\log\left(\frac{p(\bs Y|\bs{\theta})}{p(\bs Y)}\right)\right]+\frac{1}{4NM_2^2}\mathbb{V}\left[\frac{\mathbb{V}[p(\bs Y|\bs{\vartheta},\bs{\varphi})|\bs Y]}{p^2(\bs Y)}\right]+\frac{1}{4NM_1^2}\mathbb{V}\left[\frac{\mathbb{V}[p(\bs Y|\bs{\theta},\bs{\varphi})|\bs Y,\bs{\theta}]}{p^2(\bs Y|\bs{\theta})}\right]\\\nonumber
{}&+\frac{1}{NM_2}Cov\left[\log\left(\frac{p(\bs Y|\bs{\theta})}{p(\bs Y)}\right),\frac{\mathbb{V}[p(\bs Y|\bs{\vartheta},\bs{\varphi})|\bs Y]}{p^2(\bs Y)}\right]\\\nonumber
{}&-\frac{1}{NM_1}Cov\left[\log\left(\frac{p(\bs Y|\bs{\theta})}{p(\bs Y)}\right),\frac{\mathbb{V}[p(\bs Y|\bs{\theta},\bs{\varphi})|\bs Y,\bs{\theta}]}{p^2(\bs Y|\bs{\theta})}\right]\\
{}&-\frac{1}{4NM_1M_2}Cov\left[\frac{\mathbb{V}[p(\bs Y|\bs{\vartheta},\bs{\varphi})|\bs Y]}{p^2(\bs Y)},\frac{\mathbb{V}[p(\bs Y|\bs{\theta},\bs{\varphi})|\bs Y,\bs{\theta}]}{p^2(\bs Y|\bs{\theta})}\right]+\mathcal{O}_{\bb{P}}\left(\frac{1}{NM_1^2}\right)+\mathcal{O}_{\bb{P}}\left(\frac{1}{NM_2^2}\right).
\end{align}
This yields
\begin{equation}\label{d3}
  \frac{1}{N}\mathbb{V}\left[\log\left(\frac{p(\bs Y|\bs{\theta})}{p(\bs Y)}\right)\right]+\mathcal{O}_{\bb{P}}\left(\frac{1}{NM_1}\right)+\mathcal{O}_{\bb{P}}\left(\frac{1}{NM_2}\right).
\end{equation}
From \eqref{ev}, we obtain
\[\frac{1}{N}\mathbb{E}\left[\mathbb{V}\left[\log(\hat{p}_{M_1}(\bs Y|\bs{\theta}))-\log(\hat{p}_{M_2}(\bs Y))|\bs Y, \bs \theta\right]\right]=\frac{1}{N}\mathbb{E}[\mathbb{V}[\log(\hat{p}_{M_1}(\bs Y|\bs{\theta}))|\bs Y,\bs{\theta}]]+\frac{1}{N}\bb{E}[\mathbb{V}[\log(\hat{p}_{M_2}(\bs Y))|\bs Y]]\]
as the inner samples are independent.
For the first term, by the first-order Taylor expansion \eqref{tayl}, we obtain
\begin{align}
\mathbb{V}{}&[\log(\hat{p}_{M_1}(\bs Y|\bs{\theta}))|\bs Y,\bs{\theta}]\\
{}&=\mathbb{V}\left[\log(p(\bs Y|\bs{\theta}))+\frac{1}{p(\bs Y|\bs{\theta})}\left(\frac{1}{M_1}\sum_{m=1}^{M_1}p(\bs Y|\bs{\theta},\bs{\varphi}^{(m)})-p(\bs Y|\bs{\theta})\right)\Bigg|\bs Y,\bs{\theta}\right]+\mathcal{O}_{\bb{P}}\left(\frac{1}{M_1^2}\right)\\
{}&=\frac{1}{M_1}\frac{\mathbb{V}[p(\bs Y|\bs{\theta},\bs{\varphi})|\bs Y,\bs{\theta}]}{p^2(\bs Y|\bs{\theta})}+\mathcal{O}_{\bb{P}}\left(\frac{1}{M_1^2}\right).
\end{align}
Similarly, for the second term, we obtain
\begin{equation}
\mathbb{V}[\log(\hat{p}_{M_2}(\bs Y))|\bs Y]=\frac{1}{M_2}\frac{\mathbb{V}[p(\bs Y|\bs{\vartheta},\bs{\varphi})|\bs Y]}{p^2(\bs Y)}+\mathcal{O}_{\bb{P}}\left(\frac{1}{M_2^2}\right),
\end{equation}
resulting in
\[\eqref{ev}=\frac{1}{NM_1}\mathbb{E}\left[\frac{\mathbb{V}[p(\bs Y|\bs{\theta},\bs{\varphi})|\bs Y,\bs{\theta}]}{p^2(\bs Y|\bs{\theta})}\right]+\frac{1}{NM_2}\mathbb{E}\left[\frac{\mathbb{V}[p(\bs Y|\bs{\vartheta},\bs{\varphi})|\bs Y]}{p^2(\bs Y)}\right]+\mathcal{O}_{\bb{P}}\left(\frac{1}{NM_1^2}\right)+\mathcal{O}_{\bb{P}}\left(\frac{1}{NM_2^2}\right),\]
completing the derivation.

\begin{rmk}[Covariance terms in the statistical error approximation]
Applying the inequality
\begin{equation}
    \bb{V}[A+B]\leq 2\bb{V}[A]+2\bb{V}[B]
\end{equation}
demonstrates that the covariance terms in \eqref{eq:cov.terms} are of a similar magnitude as the remaining variance terms. Moreover, the covariance terms are challenging to estimate numerically and were thus neglected in the optimization of the number of samples (see \cite{Bec18}). Figure~\ref{fig:ex1.evt} (Panel A) demonstrates that this simplification did not significantly impact the consistency of the DLMC2IS estimator.
\end{rmk}
\section{Optimal setting for the double-loop Monte Carlo estimator with two inner loops}\label{ap:optimal.setting}
We solve the minimization problem stated in \eqref{eq:subproblem} by introducing a Lagrange function, which we can take the derivative of with respect to $N$, $M_1$, $M_2$, the error splitting parameter $\kappa$, and the Lagrange multipliers $\lambda$ and $\mu$ and set the resulting equations to 0\footnote{Denoted by $\stackrel{!}{=}0$.}
\begin{equation}
\mathcal{L}(N,M_1,M_2,\kappa,\lambda,\mu)=N(M_1+M_2)-\lambda\left(\frac{C_1}{M_1}+\frac{C_2}{M_2}- (1-\kappa)TOL\right)-\mu\left(\frac{1}{N}\left(D_3+\frac{D_1}{M_1}+\frac{D_2}{M_2}\right)-\left(\frac{\kappa TOL}{C_{\alpha}}\right)^{2}\right),
\end{equation}
with derivatives
\begin{equation}\label{ddn}
\frac{\partial \mathcal{L}}{\partial N}=M_1+M_2+\frac{\mu}{N^2}\left(D_3+\frac{D_1}{M_1}+\frac{D_2}{M_2}\right)\stackrel{!}{=}0,
\end{equation}
\begin{equation}\label{ddm1}
\frac{\partial \mathcal{L}}{\partial M_1}=N+\lambda\frac{C_1}{M_1^2}+\mu\frac{D_1}{NM_1^2}\stackrel{!}{=}0,
\end{equation}
\begin{equation}\label{ddm2}
\frac{\partial \mathcal{L}}{\partial M_2}=N+\lambda\frac{C_2}{M_2^2}+\mu\frac{D_2}{NM_2^2}\stackrel{!}{=}0,
\end{equation}
\begin{equation}\label{ddk}
\frac{\partial \mathcal{L}}{\partial \kappa}=2\mu\kappa\frac{TOL^2}{C_{\alpha}^2}-\lambda TOL\stackrel{!}{=}0,
\end{equation}
\begin{equation}\label{ddl}
\frac{\partial \mathcal{L}}{\partial \lambda}=(1-\kappa)TOL-\frac{C_1}{M_1}-\frac{C_2}{M_2}\stackrel{!}{=}0,
\end{equation}
and
\begin{equation}\label{ddm}
\frac{\partial \mathcal{L}}{\partial \mu}=\left(\frac{\kappa TOL}{C_{\alpha}}\right)^{2}-\frac{1}{N}\left(D_3+\frac{D_1}{M_1}+\frac{D_2}{M_2}\right)\stackrel{!}{=}0.
\end{equation}
From \eqref{eq:C1}--\eqref{eq:D2}, we have that $D_1=2C_1$ and that $D_2=2C_2$. Thus, it follows from \eqref{ddm1} that
\begin{equation}\label{eq:M1.kappa}
    \frac{M_1^2}{C_1}=-\frac{1}{N}\left(\lambda+\frac{2\mu}{N}\right),
\end{equation}
and from \eqref{ddm2} that
\begin{equation}
    \frac{M_2^2}{C_2}=-\frac{1}{N}\left(\lambda+\frac{2\mu}{N}\right).
\end{equation}
Thus, it follows that
\begin{equation}
    M_2=M_1\frac{\sqrt{C_2}}{\sqrt{C_1}}.
\end{equation}
Next, from \eqref{ddl}, we have that
\begin{align}
    (1-\kappa)TOL{}&=\frac{C_1}{M_1}+\frac{C_2}{M_2},\nonumber\\
    {}&=\frac{C_1}{M_1}+\frac{\sqrt{C_1}\sqrt{C_2}}{M_1},
\end{align}
and thus that
\begin{equation}\label{eq:M1}
    M_1^{\ast}=\frac{\sqrt{C_1}(\sqrt{C_1}+\sqrt{C_2})}{(1-\kappa)TOL},
\end{equation}
and also that
\begin{equation}\label{eq:M2}
    M_2^{\ast}=\frac{\sqrt{C_2}(\sqrt{C_1}+\sqrt{C_2})}{(1-\kappa)TOL}.
\end{equation}
Moreover, from \eqref{ddk}, we have that
\begin{equation}\label{eq:lambda.mu}
    \lambda=\frac{2\mu\kappa TOL}{C_{\alpha}^2}.
\end{equation}
Substituting \eqref{eq:M1} and \eqref{eq:M2} into \eqref{ddn}, we find that
\begin{align}
    M_1^{\ast}+M_2^{\ast}+\frac{2\mu}{N^2}\left(\frac{D_3}{2}+\frac{C_1}{M_1}+\frac{C_2}{M_2}\right){}&=\frac{(\sqrt{C_1}+\sqrt{C_2})^2}{(1-\kappa)TOL}+\frac{2\mu}{N^2}\left(\frac{D_3}{2}+(1-\kappa)TOL\right),\nonumber\\
    {}&=0,
\end{align}
yielding
\begin{equation}\label{eq:mu.N}
    \mu=-\frac{N^2}{2}\frac{(\sqrt{C_1}+\sqrt{C_2})^2}{(1-\kappa)TOL}\left(\frac{D_3}{2}+(1-\kappa)TOL\right)^{-1}.
\end{equation}
substituting into \eqref{eq:lambda.mu} leads to
\begin{equation}
    \lambda=-\frac{N^2\kappa}{C_{\alpha}^2}\frac{(\sqrt{C_1}+\sqrt{C_2})^2}{(1-\kappa)}\left(\frac{D_3}{2}+(1-\kappa)TOL\right)^{-1}.
\end{equation}
From \eqref{ddm}, we obtain that
\begin{equation}
    N^{\ast}=2\left(\frac{\kappa TOL}{C_{\alpha}}\right)^{-2}\left(\frac{D_3}{2}+(1-\kappa)TOL\right).
\end{equation}
Finally, substituting \eqref{eq:M1}, \eqref{eq:lambda.mu} and \eqref{eq:mu.N} into \eqref{eq:M1.kappa}, it follows that
\begin{align}
    \frac{\left(\sqrt{C_1}+\sqrt{C_2}\right)^2}{(1-\kappa)^2TOL^2}&{}=-\frac{1}{N^{\ast}}\left(\lambda+\frac{2\mu}{N^{\ast}}\right),\nonumber\\
    {}&=\frac{1}{N^{\ast}}\left(\frac{N^{\ast2}\kappa}{C_{\alpha}^2}\frac{(\sqrt{C_1}+\sqrt{C_2})^2}{(1-\kappa)} +\frac{N^{\ast}(\sqrt{C_1}+\sqrt{C_2})^2}{(1-\kappa)TOL}\right)\left(\frac{D_3}{2}+(1-\kappa)TOL\right)^{-1},\nonumber\\
    {}&=\left(\frac{N^{\ast}\kappa}{C_{\alpha}^2}\frac{(\sqrt{C_1}+\sqrt{C_2})^2}{(1-\kappa)} +\frac{(\sqrt{C_1}+\sqrt{C_2})^2}{(1-\kappa)TOL}\right)\left(\frac{D_3}{2}+(1-\kappa)TOL\right)^{-1}.
\end{align}
Cancellations lead to the following quadratic equation for $\kappa^{\ast}$:
\begin{align}
    1&{}=\left(\frac{N^{\ast}\kappa(1-\kappa) TOL^2}{C_{\alpha}^2}+(1-\kappa)TOL\right)\left(\frac{D_3}{2}+(1-\kappa)TOL\right)^{-1},\nonumber\\
    &{}=\left(\frac{2(1-\kappa)}{\kappa}\left(\frac{D_3}{2}+(1-\kappa)TOL\right)+(1-\kappa)TOL\right)\left(\frac{D_3}{2}+(1-\kappa)TOL\right)^{-1}.
\end{align}
Rearranging terms yields
\begin{align}
    \frac{D_3}{2}+(1-\kappa)TOL{}&=\frac{2(1-\kappa)}{\kappa}\left(\frac{D_3}{2}+(1-\kappa)TOL\right)+(1-\kappa)TOL,
\end{align}
and ultimately
\begin{equation}
    D_3\kappa=2D_3(1-\kappa)+4(1-\kappa)^2TOL,
\end{equation}
or, in the standard form
\begin{equation}
    \left(\frac{4TOL}{D_3}\right)\kappa^2-\left(3+\frac{8TOL}{D_3}\right)\kappa+2+\frac{4TOL}{D_3}=0
\end{equation}
with the solutions
\begin{equation}
    \kappa_{1,2}=1+\frac{3D_3}{8TOL}\pm\sqrt{\frac{D_3}{4TOL}+\frac{9D_3^2}{64TOL^2}}.
\end{equation}
From the requirement that $0<\kappa<1$, it follows that only the solution
\begin{equation}
    \kappa^{\ast}=1+\frac{3D_3}{8TOL}-\sqrt{\frac{D_3}{4TOL}+\frac{9D_3^2}{64TOL^2}}
\end{equation}
is permissible. It follows immediately that
\begin{equation}
    \frac{3D_3}{8TOL}<\sqrt{\frac{D_3}{4TOL}+\frac{9D_3^2}{64TOL^2}}
\end{equation}
for any $D_3,TOL>0$, and thus that $\kappa^{\ast}<1$. Moreover, it follows that
\begin{align}
    \frac{D_3}{4TOL}+\frac{9D_3^2}{64TOL^2}{}&<\left(1+\frac{3D_3}{8TOL}\right)^2,\nonumber\\
    {}&=1+\frac{3D_3}{4TOL}+\frac{9D_3^2}{64TOL^2}
\end{align}
for any $D_3,TOL>0$, and thus that $0<\kappa^{\ast}$. Applying L'H\^{o}pital's rule, we observe that
\begin{align}
  \lim_{TOL\rightarrow 0}\kappa^{\ast}={}&\lim_{TOL\rightarrow 0}1+\frac{3D_3}{8TOL}-\sqrt{\frac{D_3}{4TOL}+\frac{9D_3^2}{64TOL^2}}\nonumber\\
  ={}&\lim_{TOL\rightarrow 0}\frac{\left(8TOL+3D_3-\sqrt{16TOLD_3+9D_3^2}\right)}{8TOL},\nonumber\\
  ={}&\lim_{TOL\rightarrow 0}\frac{\frac{\di{}}{\di{}TOL}\left(8TOL+3D_3-\sqrt{16TOLD_3+9D_3^2}\right)}{\frac{\di{}}{\di{}TOL}8TOL},\nonumber\\
  ={}&\lim_{TOL\rightarrow 0}1-\frac{D_3}{\sqrt{16TOLD_3+9D_3^2}}=\frac{2}{3},
\end{align}
and that
\begin{equation}
  \lim_{TOL\rightarrow \infty}\kappa^{\ast} = 1.
\end{equation}

\section{Optimal setting with additional discretization bias}\label{ap:optimal.setting.bias}
We solve the minimization problem in \eqref{eq:subproblem.h}:
\begin{multline}\label{eq:lagr.h}
\mathcal{L}(N,M_1,M_2,h,\kappa,\lambda,\mu)=N(M_1+M_2)h^{-\gamma}-\lambda\left(\frac{D_3}{N}+\frac{D_1}{NM_1}+\frac{D_2}{NM_2}-\left(\frac{\kappa TOL}{C_{\alpha}}\right)^2\right)\\-\mu\left(C_3h^\eta+\frac{C_1}{M_1}+\frac{C_2}{M_2}- (1-\kappa)TOL\right),
\end{multline}
with the following derivatives, which we set to 0:
\begin{equation}\label{ddn.bias}
\frac{\partial \mathcal{L}}{\partial N}=(M_1+M_2)h^{-\gamma}+\frac{\lambda D_3}{N^2}+\frac{\lambda D_1}{N^2M_1}+\frac{\lambda D_2}{N^2M_2}\stackrel{!}{=}0,
\end{equation}
\begin{equation}\label{ddm1.bias}
\frac{\partial \mathcal{L}}{\partial M_1}=Nh^{-\gamma}+\frac{\lambda D_1}{NM_1^2}+\frac{\mu C_1}{M_1^2}\stackrel{!}{=}0,
\end{equation}
\begin{equation}\label{ddm2.bias}
\frac{\partial \mathcal{L}}{\partial M_2}=Nh^{-\gamma}+\frac{\lambda D_2}{NM_2^2}+\frac{\mu C_2}{M_2^2}\stackrel{!}{=}0,
\end{equation}
\begin{equation}\label{ddh.bias}
  \frac{\partial \mathcal{L}}{\partial h}=-\gamma N h^{-\gamma-1}(M_1+M_2)-\mu C_3\eta h^{\eta-1}\stackrel{!}{=}0,
\end{equation}
\begin{equation}\label{ddk.bias}
\frac{\partial \mathcal{L}}{\partial \kappa}=2\lambda\kappa \left(\frac{TOL}{C_{\alpha}}\right)^2-\mu TOL\stackrel{!}{=}0,
\end{equation}
\begin{equation}\label{ddl.bias}
\frac{\partial \mathcal{L}}{\partial \lambda}=\left(\frac{\kappa TOL}{C_{\alpha}}\right)^2-\frac{D_3}{N}-\frac{D_1}{NM_1}-\frac{D_2}{NM_2}\stackrel{!}{=}0,
\end{equation}
and
\begin{equation}\label{ddmu.bias}
\frac{\partial \mathcal{L}}{\partial \mu}=(1-\kappa)TOL-C_3h^{\eta}-\frac{C_1}{M_1}-\frac{C_2}{M_2}\stackrel{!}{=}0.
\end{equation}
The idea is to express every parameter of the Lagrangian \eqref{eq:lagr.h} as a function of the splitting parameter $\kappa$ and solve the remaining quadratic equation for $\kappa$.
Subtracting \eqref{ddm1.bias} from \eqref{ddm2.bias} results in
\begin{equation}\label{eq:M_2}
M_2=\sqrt{\frac{C_2}{C_1}}M_1,
\end{equation}
which we use in \eqref{ddl.bias} to obtain
\begin{align}\label{eq:M_1}
  M_1={}&\frac{\left(D_1+\sqrt{\frac{C_1}{C_2}}D_2\right)}{N\left(\frac{\kappa TOL}{C_{\alpha}}\right)^2-D_3},\nonumber\\
  ={}&\frac{2C_1\left(1+\sqrt{\frac{C_2}{C_1}}\right)}{N\left(\frac{\kappa TOL}{C_{\alpha}}\right)^2-D_3}.
\end{align}
The last equation follows from the definitions of $C_1$, $C_2$, $D_1$, and $D_2$.
From \eqref{ddmu.bias} we obtain
\begin{align}\label{eq:h.bias}
  h={}&\left(\frac{(1-\kappa)TOL-\frac{C_1}{M_1}\left(1+\sqrt{\frac{C_2}{C_1}}\right)}{C_3}\right)^{\frac{1}{\eta}},\nonumber\\
  ={}&\left(\frac{(1-\kappa)TOL-\frac{1}{2}\left(N\left(\frac{\kappa TOL}{C_{\alpha}}\right)^2-D_3\right)}{C_3}\right)^{\frac{1}{\eta}}.
\end{align}
From \eqref{ddh.bias}, we obtain
\begin{equation}
  \mu = -\frac{\gamma N M_1 \left(1+\sqrt{\frac{C_2}{C_1}}\right)}{\eta C_3 h^{\gamma+\eta}}.
\end{equation}
Inserting this into \eqref{ddk.bias} results in
\begin{equation}\label{eq:lambda}
  \lambda = -\frac{\gamma N M_1 \left(1+\sqrt{\frac{C_2}{C_1}}\right)C_{\alpha}^2}{\eta C_3 h^{\gamma+\eta}2\kappa TOL}.
\end{equation}
Inserting \eqref{eq:M_2} to \eqref{eq:lambda} into \eqref{ddn.bias}, after some reordering, yields
\begin{equation}
  N^{\ast}=\frac{C_{\alpha}^2}{\kappa^2TOL}\left(\frac{D_3}{TOL}+2\left(1-\kappa\left(1+\frac{\gamma}{2\eta}\right)\right)\right).
\end{equation}
This result is the same optimal number of outer samples $N^{\ast}$ as in the case for only one inner loop. We can insert this into \eqref{eq:M_1} to obtain
\begin{equation}
  M_1^{\ast}=\frac{C_1\left(1+\sqrt{\frac{C_2}{C_1}}\right)}{\left(1-\kappa\left(1+\frac{\gamma}{2\eta}\right)\right)TOL}
\end{equation}
 and
 \begin{equation}
   M_2^{\ast}=\frac{C_2\left(1+\sqrt{\frac{C_1}{C_2}}\right)}{\left(1-\kappa\left(1+\frac{\gamma}{2\eta}\right)\right)TOL}.
 \end{equation}
Equation~\eqref{eq:h.bias} results in
\begin{equation}
  h^{\ast}=\left(\frac{\gamma\kappa TOL}{2\eta C_3}\right)^{\frac{1}{\eta}}.
\end{equation}
This quantity also remains unchanged compared to the case with only one inner loop.

Finally, \eqref{ddm1.bias} provides the following quadratic equation:
\begin{equation}
  \left[\frac{1}{D_3}\left(1+\frac{\gamma}{2\eta}\right)^2TOL\right]\kappa^{\ast 2}-\left[\frac{1}{4}+\left(\frac{1}{2}+\frac{2}{D_3}TOL\right)\left(1+\frac{\gamma}{2\eta}\right)\right]\kappa^{\ast}+\left[\frac{1}{2}+\frac{1}{D_3}TOL\right],
\end{equation}
coinciding with the equation in the case with only one inner loop.

\section{Derivation of the order of the additional terms when accounting for nuisance uncertainty}\label{ap:order.maximizer}
We demonstrate the order of $\nabla_{\bs{\theta}}\bs{z}(\bs{\hat{\theta}})$, $\bs{z}(\bs{\hat{\theta}})=(\bs{\hat{\theta}},\bs{\hat{\phi}}(\bs{\hat{\theta}}))$.
We let
\begin{align}
	S(\bs{\theta},\bs{\phi})={}&\frac{1}{2}N_e(\bs g(\bs{\theta}_t,\bs{\phi}_t)-\bs g(\bs{\theta},\bs{\phi}))\cdot\bs{\Sigma}_{\bs{\varepsilon}}^{-1}(\bs g(\bs{\theta}_t,\bs{\phi}_t)-\bs g(\bs{\theta},\bs{\phi}))\nonumber\\
  {}&+\sum_{i=1}^{N_e}\bs\epsilon_i\cdot\bs{\Sigma}_{\bs{\varepsilon}}^{-1}(\bs g(\bs{\theta}_t,\bs{\phi}_t)-\bs g(\bs{\theta},\bs{\phi}))-\log(\pi(\bs{\phi}|\bs{\theta})),
\end{align}
with the gradient
\begin{equation}
	\nabla_{\bs{\phi}} S(\bs{\theta},\bs{\phi})=-N_e\nabla_{\bs{\phi}}\bs{g}(\bs{\theta},\bs{\phi})^{\trans}\bs{\Sigma}_{\bs{\varepsilon}}^{-1}(\bs g(\bs{\theta}_t,\bs{\phi}_t)-\bs g(\bs{\theta},\bs{\phi}))-\sum_{i=1}^{N_e}\bs\epsilon_i\cdot\bs{\Sigma}_{\bs{\varepsilon}}^{-1}\nabla_{\bs{\phi}}\bs{g}(\bs{\theta},\bs{\phi})-\nabla_{\bs\phi}\log(\pi(\bs{\phi}|\bs{\theta})).
\end{equation}
The first-order approximation of $\nabla_{\bs{\phi}} S(\bs{\theta},\bs{\phi})$ around $\bs{\phi}_t$, given by
\begin{equation}
	\nabla_{\bs{\phi}} S(\bs{\theta},\bs{\phi})\approx \nabla_{\bs{\phi}} S(\bs{\theta},\bs{\phi}_t)+\nabla_{\bs{\phi}}\nabla_{\bs{\phi}} S(\bs{\theta},\bs{\phi}_t)(\bs{\phi}-\bs{\phi}_t),
\end{equation}
yields
\begin{align}
	\nabla_{\bs{\phi}} S(\bs{\theta},\bs{\phi})\approx{}& \left(-N_e\nabla_{\bs{\phi}}\bs{g}(\bs{\theta},\bs{\phi}_t)^{\trans}\bs{\Sigma}_{\bs{\varepsilon}}^{-1}(\bs g(\bs{\theta}_t,\bs{\phi}_t)-\bs g(\bs{\theta},\bs{\phi}_t))-\sum_{i=1}^{N_e}\bs\epsilon_i\cdot\bs{\Sigma}_{\bs{\varepsilon}}^{-1}\nabla_{\bs{\phi}}\bs{g}(\bs{\theta},\bs{\phi}_t)-\nabla_{\bs{\phi}}\log(\pi(\bs{\phi}_t|\bs{\theta}))\right)\nonumber\\
	{}&+\Bigg(-N_e\nabla_{\bs{\phi}}\nabla_{\bs{\phi}}\bs{g}(\bs{\theta},\bs{\phi}_t)^{\trans}\bs{\Sigma}_{\bs{\varepsilon}}^{-1}(\bs g(\bs{\theta}_t,\bs{\phi}_t)-\bs{g}(\bs{\theta},\bs{\phi}_t))+N_e\nabla_{\bs{\phi}}\bs{g}(\bs{\theta},\bs{\phi}_t)^{\trans}\bs{\Sigma}_{\bs{\varepsilon}}^{-1}\nabla_{\bs{\phi}}\bs{g}(\bs{\theta},\bs{\phi}_t)\nonumber\\
	{}&-\sum_{i=1}^{N_e}\bs\epsilon_i\cdot\bs{\Sigma}_{\bs{\varepsilon}}^{-1}\nabla_{\bs{\phi}}\nabla_{\bs{\phi}}\bs{g}(\bs{\theta},\bs{\phi}_t)
	-\nabla_{\bs\phi}\nabla_{\bs\phi}\log(\pi(\bs{\phi}_t|\bs{\theta}))\Bigg)(\bs{\phi}-\bs{\phi}_t).
\end{align}

Applying Newton's method and $\nabla_{\bs{\phi}} S(\bs{\theta},\bs{\hat{\phi}}(\bs{\theta}))=\bs{0}$ implies that
$\nabla_{\bs{\phi}} S(\bs{\theta},\bs{\phi}_t)+\nabla_{\bs{\phi}}\nabla_{\bs{\phi}} S(\bs{\theta},\bs{\phi}_t)(\hat{\bs{\phi}}(\bs{\theta})-\bs{\phi}_t)\approx\bs{0}$; therefore
\begin{equation}
	\bs{\hat{\phi}}(\bs{\theta})\approx\bs{\phi}_t-\nabla_{\bs{\phi}}\nabla_{\bs{\phi}} S(\bs{\theta},\bs{\phi}_t)^{-1}\nabla_{\bs{\phi}} S(\bs{\theta},\bs{\phi}_t),
\end{equation}
resulting in
\begin{align}\label{eq:hat.phi.approx}
	\bs{\hat{\phi}}(\bs{\theta})\approx{}{}&\bs{\phi}_t\nonumber\\
	{}&-\Bigg(-N_e\nabla_{\bs{\phi}}\nabla_{\bs{\phi}}\bs{g}(\bs{\theta},\bs{\phi}_t)^{\trans}\bs{\Sigma}_{\bs{\varepsilon}}^{-1}(\bs g(\bs{\theta}_t,\bs{\phi}_t)-\bs{g}(\bs{\theta},\bs{\phi}_t))+N_e\nabla_{\bs{\phi}}\bs{g}(\bs{\theta},\bs{\phi}_t)^{\trans}\bs{\Sigma}_{\bs{\varepsilon}}^{-1}\nabla_{\bs{\phi}}\bs{g}(\bs{\theta},\bs{\phi}_t)\nonumber\\
	{}&-\sum_{i=1}^{N_e}\bs\epsilon_i\cdot\bs{\Sigma}_{\bs{\varepsilon}}^{-1}\nabla_{\bs{\phi}}\nabla_{\bs{\phi}}\bs{g}(\bs{\theta},\bs{\phi}_t)
	-\nabla_{\bs\phi}\nabla_{\bs\phi}\log(\pi(\bs{\phi}_t|\bs{\theta}))\Bigg)^{-1}\nonumber\\
	{}&\times\left(-N_e\nabla_{\bs{\phi}}\bs{g}(\bs{\theta},\bs{\phi}_t)^{\trans}\bs{\Sigma}_{\bs{\varepsilon}}^{-1}(\bs g(\bs{\theta}_t,\bs{\phi}_t)-\bs g(\bs{\theta},\bs{\phi}_t))-\sum_{i=1}^{N_e}\bs\epsilon_i\cdot\bs{\Sigma}_{\bs{\varepsilon}}^{-1}\nabla_{\bs{\phi}}\bs{g}(\bs{\theta},\bs{\phi}_t)-\nabla_{\bs{\phi}}\log(\pi(\bs{\phi}_t|\bs{\theta}))\right).
\end{align}
Both terms in the product have leading order $N_e$; thus, assuming that Newton's method converges, the Jacobian
\begin{equation}
\nabla_{\bs{\theta}}\bs{z}(\bs{\hat{\theta}})=\begin{pmatrix}I_{d_{\theta}\times d_{\theta}}\\\nabla_{\bs{\theta}}\bs{\hat{\phi}}(\bs{\hat{\theta}})\end{pmatrix}
\end{equation}
has leading order terms that are constant in $N_e$ by the quotient rule. The term $k(\bs{\theta})$ is of order $\cl{O}_{\bb{P}}\left(\log(N_e)\right)$ and the term $\ell(\bs{\theta})$ is constant in $N_e$. Evaluating $\bs{\hat{\phi}}(\bs{\theta})$ at $\bs{\theta}_t$ rather than at $\bs{\hat{\theta}}$ reduces the leading order of the numerator in \eqref{eq:hat.phi.approx} to $\sqrt{N_e}$ (see \cite{Lon13}). This would suggest approximating $(\bs{\hat{\theta}},\bs{\hat{\phi}})$ by $(\bs{\theta}_t,\bs{\phi}_t)$, however, numerical experiments show that the influence from the nuisance uncertainty would not be captured accurately by such an approximation.

%-------------------------------------------------------------------------------%

\end{document}